\documentclass[11pt]{amsart}
\usepackage{amssymb}
\usepackage{amsmath}
\usepackage{amsthm}
\usepackage{epsfig}
\usepackage{mathrsfs}
\usepackage{graphicx}
\usepackage{ textcomp }
\usepackage{xcolor}
\usepackage{hyperref,comment}
\usepackage{bbold}
\usepackage{tikz}
\usepackage{amssymb}
\usepackage{fullpage}

\numberwithin{equation}{section}

\def\RR{\mathbb{R}}
\def\NN{\mathbb{N}}

\def\e{\epsilon}
\def\G{\mathcal G}

\def\div{\mathrm {div}\,}
\def\supp{\mathrm {supp}\,}

\def\to{\rightarrow}

\def\brho{\bar \rho}

\def\vphi{\varphi}
\def\pa{\partial}
\def\na{\nabla}
\def\eps{\varepsilon}

\def\ds{\displaystyle}

\newtheorem{theorem}{Theorem}[section]
\newtheorem{lemma}[theorem]{Lemma}
\newtheorem{proposition}[theorem]{Proposition}
\newtheorem{definition}[theorem]{Definition}
\newtheorem{assumption}[theorem]{Assumption}

\newtheorem{remark}[theorem]{Remark}
\newtheorem{rk&ex}[theorem]{Remarks \& Examples}
\newtheorem{corollary}[theorem]{Corollary}

\begin{document}
%\title{Incompressible  limit of  porous medium equation with a sign-changing growth term}
\title{A Hele-Shaw limit without monotonicity}
\author[N. Guillen]{Nestor Guillen}
\address{Department of Mathematics, Texas State University, San Marcos, TX}
\email{nestor@txstate.edu}

\author[I. Kim]{Inwon Kim}
\address{Department of Mathematics, UCLA,  Los Angeles, CA} 
\email{ikim@math.ucla.edu}

\author[A.Mellet]{Antoine Mellet}
\address{Department of Mathematics, University of  Maryland, College Park, MD}
\email{mellet@umd.edu}

\begin{abstract}
We study the incompressible limit of the porous medium equation with a right hand side representing either a source or a sink term, and an injection boundary condition. This  model can be seen as a simplified description of non-monotone motions in tumor growth and crowd motion, generalizing the congestion-only motions studied in recent literature (\cite{AKY}, \cite{PQV}, \cite{KP}, \cite{MPQ}). We characterize the limit density, which solves a free boundary problem of Hele-Shaw type in terms of the limit pressure.  The novel feature of our result lies in the characterization of the limit pressure, which solves an obstacle problem at each time in the evolution.
\end{abstract}

\maketitle

\section{Introduction}
The porous media equation is a nonlinear evolution equation which is commonly used to model many natural phenomena involving diffusion or heat propagation. In its simplest form, it consists of a continuity equation with a flux given by Darcy's law:
\begin{equation}\label{pressure} 
\pa_t  \rho - \div(\rho \nabla p) = 0, \quad p = \frac{m}{m-1} \rho^{m-1}, \qquad m>1.
 \end{equation}
The exponent $m>1$ describes the anti-crowd tendency of the density motion, where the diffusion is larger at higher density (\cite{BGHP}, \cite{BH}, \cite{M}, \cite{W}). Due to the degeneracy of the diffusion at lower densities, it is well-known that the density stays compactly supported if initially so (see for example \cite[Chapter 1]{Vaz2007}). 

\medskip

In this paper, we consider the porous media equation with a source term:
\begin{equation}\label{eq:0}
\pa_t \rho- \div(\rho\na p) = \lambda \rho \qquad  \mbox{ in }  \Omega\times\RR^+,
\end{equation}
set in the (exterior) domain $\Omega = \RR^n\setminus K$, where $K$ is a bounded subset of $\RR^n$ with smooth boundary and supplemented with the ``injection" boundary condition,
\begin{equation}\label{eq:bc}
 \rho(x,t)=f(x,t)^{\frac{1}{m-1}} > 0\quad \mbox{ on } \pa K.
 \end{equation}
 as well as the initial data $\rho(t=0)=\rho_m^0\geq 0$. 
 Importantly we will assume that the function $\lambda(x,t)$ is bounded but can take both positive or negative values.
When $\lambda <0$, the term $\lambda \rho$ is an absorption term which is competing with the injection at $\pa K$.
Assumptions on the  initial and boundary data are given in Section \ref{sec:ass}. 
Classically, \eqref{eq:0} can also be written as the following equation for the pressure  $p(x,t)$:
 \begin{equation}\label{eq:pressure}
\pa_t p = (m-1)p (\Delta p +\lambda) + |\na p|^2.
\end{equation}
  
 \medskip
 
Our interest is with the {\it incompressible limit} of this equation, that is the limit $m\to\infty$.  
Heuristically speaking, if $(\rho_m,p_m)$ denotes a sequence of solution of \eqref{eq:0}, then -- provided there is an actual limit in a good enough sense -- the limits $\rho_\infty$ and $p_\infty$ should satisfy
\begin{equation}\label{eq:limit1}
\pa_t \rho_\infty - \div(\rho_\infty \nabla p_\infty) = \lambda \rho_\infty\quad  \mbox{ in } \Omega\times\RR^+, \qquad p_\infty = f\quad \mbox{ on } \pa K, \qquad \rho_{\infty}(\cdot,0) = \rho^0,
\end{equation}
and  taking the limit in the relation $p_m = \tfrac{m}{m-1}\rho_m^{m-1}$, we may guess that in the limit $\rho_\infty$ and $p_\infty$ are connected by what is known as the Hele-Shaw graph
\begin{equation}\label{eq:Pgraph} 
p_\infty \in P_\infty (\rho_\infty) : = \begin{cases}
0 & \mbox{ if }  0\leq \rho_\infty<1 \\
[0,\infty) & \mbox{ if } \rho_\infty=1\\
\infty & \mbox{ if } \rho_\infty>1.
\end{cases}
\end{equation}
In particular, the pressure can be viewed as a Lagrange multiplier for the constraint $\rho_\infty \leq 1$ (\cite{MRS}). 
In our framework, as in many of the related works discussed below,  a priori estimates (under appropriate assumptions on $f$ and on the initial data) will allow us to make the derivations of \eqref{eq:limit1}-\eqref{eq:Pgraph} rigorous.

\medskip

Equations \eqref{eq:limit1}-\eqref{eq:Pgraph} fully characterize the evolution of $\rho_\infty$ (see the uniqueness result, Proposition \ref{prop:uniqueness}). However, one would like to give a more 
geometrical description of the evolution of $\rho_\infty$, and in particular of the evolution of the ``saturated region''
\begin{align*}
  \Sigma(t):= \{ \rho_\infty(t) = 1\}.
\end{align*}

Classically, such a description is provided by a Hele-Shaw type free boundary problem.  Indeed, formally at least, we can pass to the limit in \eqref{eq:pressure} to get
the so-called {\it complementarity condition}:
%and assuming $\pa_t p_m$ or $|\na p_m|^2$ do not become too large as $m\to\infty$ then $p_\infty$ would also solve
\begin{align}\label{eq:limit2}
  p_\infty (\Delta p_\infty + \lambda) = 0 \textnormal{ in } \Omega\times\RR^+
\end{align}	
which implies that $p_\infty(\cdot,t)$ solves $-\Delta p_\infty = \lambda$ in the set $\{p_\infty(\cdot,t)>0\}$,
and equation \eqref{eq:limit1} implies (in a weak form) that the normal velocity of the interface $\pa\Sigma(t)$ is proportional to  $|\na p_\infty|$.
However, the derivation of \eqref{eq:limit2} is less straightforward than that of \eqref{eq:limit1}-\eqref{eq:Pgraph} in  general (see for instance \cite{DP}), and it is not obvious that we should always have $\{p_\infty(\cdot,t)>0\}= \Sigma(t)$.

\medskip

Incompressible limits were first studied for equation \eqref{pressure}  (that is when $\lambda = 0$).
%This analysis has been made rigorous mostly in the case where there is no reaction term, that is when $\lambda = 0$. 
In the absence of $K$, there are classical works starting by B\'enilan and Crandall \cite{BC1981}, followed by results with more general initial data by Caffarelli and Friedman \cite{CF1987} and numerical studies describing the shape of the limit by Elliot et al \cite{EHKO1986} (also see \cite{GQ03} for its rigorous justification).
In \cite{BMM2009} a similar weak formulation for the Hele-Shaw problem (still without right hand side) is derived as a "mesa" limit from the Stefan problem.
The last decade has seen significant advances in the study of these asymptotics when the right hand side is monotone increasing in $\rho$ -- corresponding to the case $\lambda>0$ in our framework. The convergence as $m\to\infty$ and characterization of the limit  as a Hele-Shaw type flow has been achieved for models of congested crowd motion \cite{AKY,BKP} and of tumor growth \cite{PQV,MPQ,KP}. It is important to note that monotonicity properties are present in the systems studied in these papers and are essential for proving that $\{p_\infty(\cdot,t)>0\}= \Sigma(t)$. For instance, the monotonicity of the density was a key feature in characterizing the limiting problem in \cite{MPQ,KP}. In \cite{AKY,BKP} which features a drift field, the monotonicity  of $\rho$ along the streamline was crucial to characterize the limiting problem in terms of viscosity solutions. 

\medskip

In our work, the function $\lambda(x,t)$ is not necessarily positive so that one no longer expects $\rho_\infty$ to be monotone in time thus complicating the analysis.
Moreover, a Hele-Shaw type problem with a single phase is typically monotone in time, suggesting that the lack of monotonicity should be reflected by having some modification of the one-phase Hele-Shaw model in the limit. 
%Determining this limit and in particular describing the dynamics of $\Sigma(t)$ when $\lambda$ might not have a definite sign is the main objective of this paper.
One of the main contribution of this paper is to identify the pressure $p_\infty(\cdot,t) $  for all time $t>0$ by showing that it solves an obstacle problem in the set $\Sigma(t)$ and might thus be such that $\{p_\infty(\cdot,t)>0\}\subsetneq\Sigma(t)$ (see Theorem \ref{prop:obs}), causing the saturated set $\Sigma(t)$ to shrink. 
Though our result appears to be new, its proof is relatively simple and can be generalized to problems with nonlinear source terms. As an illustration of this latter point, we apply these ideas to a tumor growth model which involves nonlinear terms (see Appendix \ref{app:tumor}). Even in the monotone cases mentioned above, our result provides a new approach to the derivation of the complementarity condition \eqref{eq:limit2}.

\medskip

\medskip

Equation \eqref{eq:0} is simple but it allows us to study a very general and important behavior. 
Indeed, the monotonicity in the aforementioned works is characteristic of systems with only congestive effects. 
However, it is clear that de-congestion effects are important for applications. In \cite{PQV}, a model for tumor growth which takes into account the evolution of the density of nutrients is introduced and studied. In that case, the tumor cells decrease their density in the event of insufficient nutrient, which yields to "de-congestion" or recession of the tumor cells. The consequent lack of monotonicity significantly complicates the analysis:  The derivation of the complementarity condition was only achieved recently \cite{DP} and the geometric description of the tumor growth still remains to be understood. 
Similarly, the study of congested crowd motion that involve de-congestion phenomena is of great interest (see \cite{MRS},\cite{S_survey}).

Our interest in studying the toy problem \eqref{eq:0} is thus to better understand such behavior. 
By allowing $\lambda$ to take both positive and negative value, we generate a motion that consists of both congestion and de-congestion. 
The presence of  a fixed boundary condition on $\pa K$ is by no mean necessary for our analysis  (there is no $K$ in the tumor growth model studied in Appendix \ref{app:tumor}), but such injection boundary conditions are a classical feature of Hele-Shaw problems. 
In the context of crowd motion, our model describes a congested crowd coming out of the door ($\partial K$) to the outdoors ($\mathbb{R}^n\setminus K$). In the context of the classical Hele-Shaw flow (with $\lambda=0$) the boundary condition describe the injection of the fluids.

\medskip

In our setting it seems natural to expect that $p_{\infty}$, which acts against congestion, may vanish even when the density is fully saturated. Indeed we will see that   when $\lambda$ is not necessarily positive the support of the pressure $p_\infty(t)$ may be a strict subset of $\Sigma(t)$.  In general, $p_\infty(t)$ must be found by solving an obstacle problem in $\Sigma(t)$.
As a result, while $\Sigma(t)$ will expand according to a Hele-Shaw type  law when $|\na p_\infty|>0$ along $\pa \Sigma(t)$, it might recede when $|\na p_\infty|=0$. 
Formally, the motion law of $\Sigma(t)$ can be written as
\begin{equation}\label{eq:fb}
 |\na p_\infty| = (1-\rho^E) V \hbox{ on } \partial\Sigma(t),
\end{equation}
where $V$ denotes the outer normal velocity of $\pa \Sigma(t)$ and $\rho^E$ is the trace of the``external density", namely the trace of $\rho_\infty$ on $\pa \Sigma(t)$ from $\{\rho_\infty<1\}$ (this is well defined if $\pa \Sigma(t)$ smooth since $\rho_\infty$ is in $\mathrm {BV}_{loc}(\RR^n\setminus K)$).

\medskip

The velocity law \eqref{eq:fb} can be formally justified from the weak equation \eqref{eq:limit1} as follows (where  $\nu$ denotes the  inward normal unit vector on $\pa K$):
\begin{align*}
\int_{\partial K } \rho \nabla p\cdot \nu \,  dS  + \int_\Omega \lambda \rho \, dx = \frac{d}{dt} \int _{\Omega} \rho \, dx &= \frac{d}{dt}\left[\int_{\Sigma(t)} \rho \, dx + \int_{\Omega\setminus \Sigma(t)} \rho^E\,  dx\right] \\ 
&= \int_{\Sigma(t)} \pa_t  \rho \, dx + \int_{\Omega\setminus \Sigma(t)} \pa_t \rho^E\,  dx +  \int_{\partial\Sigma(t)} V(1-\rho^E) \, dS \\ 
&= \int_{\Sigma(t)} \div(\rho\na p) +\lambda \rho \, dx + \int_{\Omega\setminus \Sigma(t)} \lambda \rho ^E\,  dx +  \int_{\partial\Sigma(t)} V(1-\rho^E) \, dS \\ 
&=\int_{\pa K} \rho \na p\cdot \nu \, dS + \int_{\pa \Sigma(t)} \rho \na p\cdot \nu \, dS
+\int_{\Omega} \lambda \rho \,  dx +  \int_{\partial\Sigma(t)} V(1-\rho^E) \, dS,
\end{align*}
from which we deduce (since $\rho=1$ in $\Sigma(t)$) that  $\displaystyle \int_{\partial\Sigma(t)} [\nabla p \cdot \nu+ V(1-\rho^E)] dS=0$.

\medskip

We note that our motion law is different from \cite{BKP} where the free boundary can move back and forth under the action of a force field. Here the receding and advancing behavior of the free boundary takes place via completely different mechanisms.
The motion law \eqref{eq:fb} is closer to the one obtained in \cite{KM} in the context of liquid drops sliding down on inclined plane. In this context, at the receding end of the drop, the contact angle between the liquid drop and the plane may vanish. In that moment the nature of the velocity law suddenly changes:  it is no longer dictated by the local value of the pressure, but rather by the bulk behavior of the liquid via an obstacle problem.
\medskip

Finally, we believe that our approach developed for the model problem \eqref{eq:0} is quite general and is of independent interest.
To illustrate this point we prove in Appendix \ref{app:tumor} that it can be applied to the tumor growth problem with nutrient, considered in \cite{PQV,DP}.  
%Interestingly, while \eqref{eq:limit1} has a unique weak solution, we are only able to derive the obstacle problem for the pressure variable based on the approximation \eqref{pressure}.

\medskip

Here is a brief outline of the paper. In Section 2 we collect and discuss implications of our results. In Section 3 we show convergence of the density and pressure variables. In section 4 we derive the novel characterization of the pressure via an obstacle problem. Section 5  introduces the comparison principle, as well as the uniqueness, of the limit problem, which will be used in the rest of the paper. Section 6 - 8 describes the motion law of the saturated region, starting with the measure theoretic representation in Section 6. An alternative characterization, in the flavor of viscosity solutions, is given in Sections 7-8.

\medskip

\subsection{Acknowledgements} The authors would like to acknowledge the generous support of the National Science Foundation. Inwon Kim was partially supported by National Science Foundation grant DMS-1900804 and Antoine Mellet was partially supported by National Science Foundation grant DMS-2009236.

\section{Notations and main results}
\subsection{Assumptions}\label{sec:ass}

Throughout the paper, we denote by $(\rho_m,p_m)$ the solution of the  following initial boundary value problem:
\begin{equation}\label{eq:1}
\begin{cases}
\pa_t \rho_m - \div(\rho_m\na p_m) = \lambda \rho_m  & \mbox{ in } Q, \qquad p_m = \frac{m}{m-1} \rho_m^{m-1}, \\[5pt]
\rho_m(x,t)=f(x,t)^{\frac{1}{m-1}} & \mbox{ on } \pa K \times \RR_+ \\[5pt]
\rho_m(0,x) = \rho^0_m(x) & \mbox{ in } \Omega
\end{cases}
\end{equation}
where we denote
$$\Omega  := \RR^n\setminus K, \quad Q:=\Omega\times \RR_+, \quad 
Q_T:=\Omega\times (0,T].$$

Below are the main assumptions  to be used throughout our analysis:
\begin{assumption}\label{ass:data} There is a constant  $\Lambda >0$ such that 
\item[(i)] The function $\lambda(x,t)$ satisfies
\begin{equation}\label{eq:lambdabd}
| \lambda(x,t) | \leq \Lambda \qquad \forall (x,t)\in Q,
\end{equation}
\begin{equation}\label{eq:lambdaBV}
\lambda \in BV_{loc}(\Omega\times\RR_+).
\end{equation}
\item[(ii)] The boundary data $f(x,t)$ satisfies
\begin{equation}\label{eq:condbc}
0<\Lambda^{-1}\leq f\leq \Lambda ,\qquad  |\na f| \leq C, \qquad |\pa_t f|\leq C \quad \mbox{ on } \pa K \times \RR^+.
\end{equation}
\end{assumption}

In order to write the assumptions on the initial condition $\rho^0_m$, we first introduce appropriate  barriers.
Given $0 \leq \underline R < \overline R$, we consider 
$\overline \vphi(x)$ and $\underline \vphi(x)$ solutions of 
\begin{equation}\label{eq:ophi}
-\Delta \overline \vphi = \Lambda +1 \mbox{ in } B_{\overline R}\setminus K, \quad \overline \vphi = f^\frac{m}{m-1} \mbox{ on } \pa K,\quad \overline \vphi = 0 \mbox{ on } \pa B_{\overline R}
\end{equation}
and
\begin{equation}\label{eq:uphi}
-\Delta \underline \vphi = -\Lambda \mbox{ in } B_{\underline R} \setminus K , \quad \underline \vphi = f^\frac{m}{m-1} \mbox{ on } \pa K,\quad \underline \vphi = 0  \mbox{ on  } \pa B_{\underline R},
\end{equation}
where we assume that $\underline \vphi>0$ in $B_{\underline R} \setminus K $ (if necessary we can replace $B_{\underline R}$ by a smaller set  sufficiently close to $K$).
\begin{assumption}\label{ass:init}
\item[(i)] The initial condition $\rho^0_m(x)$ satisfies 
\begin{equation}\label{eq:condin1}
\underline \vphi(x) ^{\frac 1  m} \leq \rho^0_m(x)\leq \overline \vphi^{\frac 1  m}(x) \qquad\forall x\in\Omega
\end{equation}
\begin{equation}\label{eq:condin3}
\| \Delta( {\rho^0_m}^m)  + \lambda  \rho^0_m \|_{L^1(\Omega)} + \| \na \rho^0_m\|_{L^1(\Omega)}  \leq C 
\end{equation}
\item[(ii)] The sequence $\{\rho^0_m\}_{m\geq 1}$ converges in $L^1(\Omega)$ to $\rho^0$.
\end{assumption}

\medskip

Condition \eqref{eq:condin3} may seem restrictive, but the following result shows that a wide range of initial condition $\rho^0$ can fit into this framework: 
\begin{lemma}
Let $\Sigma\supset K$ be a bounded open set with $C^2$ boundary in $\RR^n$ and let
 $\rho^0(x)$ be given by
$$
\rho^0:= \chi_{\Sigma} + \rho^E\chi_{\Sigma^C} \hbox{ in } \Omega,
$$
where $\rho^E \in C^{1,1}_c(\Omega)$ satisfies $0\leq \rho^E <1$.
Then there exists a sequence $\rho^0_m$ satisfying Assumption \ref{ass:init}.
\end{lemma}

The construction is simple, so we give it here: First, we define the pressure $p_0$ by 
$$
-\Delta p_0 = 0 \hbox{ in }\Sigma\setminus K \quad \hbox{ with } p_0=0\hbox{  on } \partial\Sigma \hbox{ and } p_0 = f \hbox{  on }K.
$$
We clearly have $p_0\geq 0$ in $K$ and $|\na p_0|\neq 0$ on $\partial\Sigma$.  We can then define
$$
\rho^{0}_m:= \max\{ p_0^{1/m}, (\rho^E - a_m)_+\}, 
$$
where $a_m$ is a nonnegative sequence such that $a_m \to 0$ and $(1-a_m)^m\to 0$ as $m\to\infty$ (for instance $a_m= (\ln m)^{-1}$).
Note that  with  this definition \eqref{eq:condin1}  holds for sufficiently large $m$. To check \eqref{eq:condin3},  first note that $p_0^{1/m}$ is in $BV$, since 
$$
\|Dp^{1/m}\|_{L^1}= \left\|\frac{p^{\frac 1 m-1}}{m} Dp\right\|_{L^1} \leq C\sup_{p \geq 1/m}  |Dp| + o(1),
$$
for sufficiently large $m$, where we have used the fact that $p$ grows at most linearly near the regular boundary $\partial\Sigma$.  Lastly, note that 
$$
(\rho^0_{m})^m= \max\{ p_0,  (\rho^E - a_m)^m\},
$$
which is a maximum of two $C^2$ functions. Moreover for large $m$ we have $\nabla  (p_0 - (\rho^E- a_m)^m)  \neq 0$ where they  coincide, since $\nabla p_0\neq  0$  due to  the regularity of $\partial\Sigma$ and $\nabla (\rho^E- a_m)^m$ uniformly vanishes as $m$ grows. This nondegeneracy   yields the regularity of  the set  $\Gamma:= \{p_0 = (\rho^E- a_m)^m\}$. Collecting the facts  we conclude  \eqref{eq:condin3}, where $\Delta {\rho^0_m}^m$ is interpreted as a measure.

\begin{comment}
These assumption can be easily satisfied if we assume that $\rho^0_m\geq 0$ is the compactly supported solution of
$$
\Delta (\rho^0_m)^{m}  = \lambda_2  \rho^0_m  \mbox { in } \Omega, \quad  \rho^0_m=f^\frac{m}{m-1} \mbox{ on } \pa K.
$$
\textcolor{black}{Is this function uniquely determined for a given compact set in $\RR^n$?}
\end{comment}

\subsection{Limit and weak formulation of the limiting problem}
By generalizing classical a priori estimates to our equation, we will first establish the convergence of $\rho_m$ and $p_m$ and prove the following result:
\begin{theorem}\label{thm:1}
Under Assumptions \ref{ass:data} and \ref{ass:init} and up to a subsequence,
the density $\rho_m$ and pressure $p_m$ solution of \eqref{eq:1}  converge strongly in 
$L^1(Q_T)$ for all $T>0$  to limits $\rho_\infty$ and $p_\infty$ which satisfy
$$ \rho_\infty , \; p_\infty\in BV(Q_T),$$
$$ \rho_\infty \in C^s([0,\infty);H^{-1}(\Omega)) \quad\forall s<1/2, \quad 
p_\infty \in L^2(0,T;H^1(\Omega)),
$$ 
$$  0\leq p_\infty(x,t)\leq C \quad \mbox{ a.e. } (x,t)\in Q,  \quad 0\leq \rho_\infty(x,t)\leq 1\quad \mbox{ a.e. } x\in\Omega, \; \forall t>0$$
and
\begin{equation}\label{eq:2}
\begin{cases}
\pa_t \rho_\infty = \Delta p_\infty + \lambda \rho_\infty & \mbox{ in } \mathcal D'(\Omega\times \RR^+), \qquad p_\infty \in P_\infty(\rho_\infty);\\
p_\infty(x,t)=f(x,t) & \mbox{ on } \pa K\times \RR^+; \\
\rho_\infty(x,0) = \rho^0(x) & \mbox{ in } \Omega,
\end{cases}
\end{equation}
where $P_\infty$ is the Hele-Shaw graph \eqref{eq:Pgraph}.
\end{theorem}
Following \cite{PQV}, we can prove the following result which shows that the result above 
 fully characterizes the function $\rho_\infty$:
\begin{proposition}\label{prop:uniqueness}
Suppose $\lambda \in L^2([0,T]; H^1(\Omega))$, then equation \eqref{eq:2} has at most one  solution $(\rho,p)\in X:=L^\infty(\Omega\times (0,T]) \times L^2(0,T;H^1(\Omega))$.
\item Furthermore, if $(\rho_1,p_1)$ and $(\rho_2,p_2)$ are respectively sub and super-solutions of  \eqref{eq:2}  in $X$ satisfying $\rho_1(\cdot,0)\leq \rho_2(\cdot,0)$ and $p_1 |_{\pa K} \leq p_2 |_{\pa K}$, then 
$ \rho_1\leq \rho_2$ in  $\Omega\times\RR^+$.
 \end{proposition}

\begin{remark}\label{rem:convergence}
This uniqueness result implies in particular that any subsequence of $(\rho_m,p_m)$ converges to the same limit, and thus the entire sequence converges to $(\rho_\infty,p_\infty)$.

\end{remark}

When $\lambda=0$, equation \eqref{eq:2} implies that the saturated region
$\Sigma(t)=\{\rho_\infty(t) =1 \}$ coincides with the set   $\{p_\infty(t)>0\}$, and $\Sigma(t)$  evolves according to the classical Hele-Shaw free boundary problem:
$$
\begin{cases}
\Delta p_\infty = 0 \mbox{ in } \Sigma(t),\qquad p_\infty = f \mbox{ on } \pa K,\qquad  p_\infty=0 \mbox{ on } \pa\Sigma(t);\\
V = |\na p_\infty| \mbox{ on } \pa \Sigma(t),
\end{cases}
$$
where $V$ denotes the outer normal velocity of the interface $\pa \Sigma(t)$ (\cite{K03}, \cite{QV99}).
This provides a simple geometric description of the evolution of  the set $ \{\rho_\infty=1\}$. 
As explained in the introduction, our goal in this paper is to provide a similar characterization when $\lambda \neq 0$.

\subsection{The pressure $p_\infty(t)$}

Our first task is to determine how the pressure $p_\infty(t)$ depends on the set $ \{\rho_\infty=1\}$. 
%We will see in the next section that the pressure does not necessarily solve a simple elliptic equation in $ \{\rho_\infty=1\}$ as above.
An important and  new feature in our framework is that that we may have  $\{p_\infty(t)>0\} \subsetneq \{\rho_\infty(t) =1 \} $. Indeed we prove that $p_\infty(t)$ is determined by solving an obstacle problem in the set $ \{\rho_\infty=1\}$.

First, we note that for all $t_0\geq 0$ we have $p_\infty \in BV(\Omega\times (t_0,T))$ and so 
 we can define the trace of the function $p_\infty$ on $\{t=t_0\}$. The interested reader might consult Giusti's book \cite[Chapter 2]{Giusti1984} for a thorough discussion on traces of BV functions (it is worth emphasizing that $p_\infty$ is of bounded variation in space \emph{and} time). 
 We denote this trace $p^+(x,t_0)$ since it is defined as a limit as $t\to t_0^+$.
 It satisfies in particular
\begin{equation}\label{eq:tracep} \frac 1 \delta \int_{t_0}^{t_0+\delta} \int_\Omega |p_\infty(x,t) - p^+(x,t_0)| \, dx \, dt \leq  \int_{t_0}^{t_0+\delta} \int_\Omega |\pa_t p_\infty| \, dx \to 0 \qquad \mbox{ as } \delta\to0.
\end{equation}
and (by Lebesgue differentiation theorem)
$ p_\infty(x,t) = p^+(x,t)$ almost everywhere. 
Since $\lambda \in BV$, we can similarly define the trace $\lambda^+(\cdot,t)$ for all $t>0$. 
We then prove:

\begin{theorem}\label{prop:obs}
Under the conditions of Theorem \ref{thm:1} and  for all $t\geq 0$, 
%$$p_{\infty}(\cdot,t)=p^*(\cdot,t)$$
$p^+(\cdot,t)$ is the unique solution of the minimization problem
\begin{equation}\label{eq:obstacle}
\underset{v\in E_t}{\rm min} \displaystyle \int_{\Omega} \frac 1 2 |\na v|^2 - \lambda^+(\cdot,t) v\, dx
\end{equation}
where $E_t$ denotes the functional space
$$ E_t = \left\{ v\in H^1(\Omega)\cap L^1(\Omega)\, ;\, v=f \mbox{ on } \pa K, \; v \geq 0 \mbox{ in } \Omega, \;    \langle v, 1-\rho_\infty(t) \rangle_{H^1,H^{-1}} =0 \right\}.$$
Equivalently, $p^+(\cdot,t)$ is the unique solution of the variational inequality 
\begin{equation}\label{eq:var}
\begin{cases}
p \in E_t \\   
\displaystyle \int_{\Omega} \na p \cdot \na (p-u) - \lambda^+(\cdot,t) (p-u) \, dx \leq 0 \qquad \forall u \in E_t.
\end{cases}
\end{equation}
\end{theorem}

\medskip

If the set $\Sigma(t)=\{\rho_\infty(\cdot,t) =1\}$ is a smooth enough subset of $\Omega$, then
\eqref{eq:obstacle} is a classical obstacle problem in $\Sigma(t)$ with Dirichlet boundary conditions $p=f$ on $\pa K$, $p=0$ on $\pa \Sigma(t)$. 
The proof of this result is surprisingly simple and quite flexible (see Section 4).
It does not require any additional a priori estimates besides the ones already used to prove Theorem \ref{thm:1}.
It can easily be adapted to more complicated models, such as the tumor growth model with nutrient, as we show in Appendix \ref{app:tumor} (see Proposition \ref{prop:obstumor}).

\medskip

By using the approach developed in \cite{MRS},  it is also possible to show that for any weak solutions of \eqref{eq:2}, the pressure $p_\infty(\cdot,t)$ satisfies, for a.e. $t>0$
$$ \int_\Omega \na p \cdot \na u - \lambda u \, dx=0, \qquad \forall u\in E_t.$$
So the pressure $p_\infty(\cdot,t)$ solves the equation $\Delta p +\lambda =0$ in the set $\{\rho_\infty(t)=1\}$ for almost every time.
As explained in the introduction, we cannot expect this to hold for all time, since either $\lambda$ or the set $\{\rho_\infty(t)=1\}$ may evolve discontinuously over time. In the event where the solution of the obstacle problem \eqref{eq:var} has its support strictly smaller that $\{\rho_\infty(t)=1\}$,  the set $\{\rho_\infty(t)=1\}$ will shrink instantaneously.
The result of \cite{MRS} does not see these instantaneous collapses (which can happen over a  large set of time, albeit one of measure zero). 
Our characterization of $p_\infty$, which holds for all time $t>0$, identifies how such collapses take place.

\medskip

When $\lambda=0$, Theorem \ref{prop:obs}  provides a simple proof of the harmonicity of  $p_\infty$ in $\{\rho_\infty(\cdot,t) =1\}$. In the general case, it implies in particular the so-called complementarity condition:
$$ p_\infty (\Delta p_\infty + \lambda) = 0 \mbox{ in } \mathcal D'(\Omega\times (0,\infty))$$
which is readily obtained by taking $u=p(1\pm \epsilon\vphi)$ in  \eqref{eq:var} with $\vphi \in \mathcal D(\Omega\times (0,\infty))$ and $\epsilon$ small enough so that $1\pm \epsilon\vphi\geq 0$. 

This complementarity condition is proved for the tumor growth model in \cite{PQV} (model without nutrient) and in \cite{DP} (model with nutrient). In both cases, the derivation relies on further estimates on the pressure (in particular the Aronson-B\'enilan estimate or some variant of it). Our result thus provides an alternative derivation of this condition that does not require any of these additional estimates.
\medskip

Given the interest for the complementarity condition in the literature, it is worth noting that it is equivalent to the obstacle problem formulation in the following sense:
\begin{proposition}\label{prop:comp}
Let $(\rho,p)\in L^\infty(0,T; L^1(\Omega)\cap L^\infty(\Omega) ) \times  L^2(0,T;H^1(\Omega))$ be a solution of \eqref{eq:2} with $p\in BV_{loc}(\Omega\times \RR_+)$.
If $p$ satisfies the complementarity condition
$$ p (\Delta p  +\lambda) = 0 \mbox{ in } \mathcal D'(\Omega\times (0,\infty))$$
then 
%$p(\cdot,t)$ is the unique solution of \eqref{eq:obstacle} for a.e. $t>0$.
for every $t>0$
the trace
$p^+(\cdot,t)$ (as defined in \eqref{eq:tracep}) is the unique solution of problem \eqref{eq:obstacle}.

\end{proposition}

Note that given a weak solution of \eqref{eq:2}, we are not able to prove directly that it 
satisfies the obstacle problem formulation of Theorem \ref{prop:obs}) or the complementarity condition, but this proposition shows that these two properties are equivalent.

\medskip

In general little is known on the boundary regularity of the set $\{\rho_\infty(\cdot,t) =1\}$, including whether its boundary has measure zero. 
 Thus for pointwise characterization of the pressure $p_\infty$, we  define the {\it support } of the measure $1-\rho_\infty$ by 
 $$
 \mathrm{Supp} \, (1-\rho_\infty(t)) := \left\{ x_0\in \Omega \,;\, \int_{B_r(x_0)}(1-\rho_\infty)(\cdot,t)\,  dx >0 \mbox{ for all } r>0
\right\} .
$$
While it may differ from the set $\{\rho_{\infty} <1\}$ by a measure zero set,  this set has the advantage of being closed by its definition. Then the solution of the obstacle problem \eqref{eq:obstacle} has the usual properties in the open set  
\begin{equation}\label{interior}
\mathcal O(t) := \Omega \setminus \mathrm{Supp} \, (1-\rho_\infty(t)) = \left\{x_0\in\Omega: \int_{B_r(x_0)} (1-\rho_\infty)(\cdot,t) = 0 \mbox{ for  some } r>0\right\}
\end{equation}
which can be seen as the ``interior" of the set $\{\rho_\infty(\cdot,t) =1\}$. 
More precisely, we have:
\begin{proposition}\label{prop:obstacle}
The function $p^*$, solution of the minimization problem \eqref{eq:obstacle},  is in $C^{1,1}_{loc}(\mathcal{O}(t))$ and satisfies
\begin{equation}\label{eq:pob}
- \Delta p^* = \lambda \chi_{\{p^*>0\}} \mbox{ in } \mathcal O(t).
 \end{equation}
\end{proposition}

\begin{comment}
 More can be said for lower dimensions due to the following regularity result, which is an adaption of similar estimates shown in \cite{MPQ} and \cite{DP} for tumor growth models:
\begin{proposition}\label{prop:L4}
For every test function $\vphi \in \mathcal D(\Omega)$, there exists a constant $C$, independent of $m$ such that
$$\int_{\frac{1}{m-1}}^T \int_\Omega |\na p_m(x,t)|^4 \vphi(x) \, dx\, dt\leq C.$$
In particular, $p_\infty(\cdot,t) $ is in $W^{1,4}_{loc}(\Omega)$ for a.e. $t>0$.
\end{proposition}
When $n\leq 3$,  the proposition shows that $p_{\infty}(\cdot,t)$ is continuous function for a.e. $t>0$. In particular $p_{\infty}$ thus solves a classical obstacle problem in the set $\mathcal O(t)$ with the dirichlet boundary condition $p_\infty(\cdot,t) = 0$ on $\pa \mathcal O(t)$. 
\end{comment}
 \medskip

\subsection{Velocity law: Measure theoretic results}
In view of Theorem \ref{prop:obs}, we can redefine $p_\infty$ a.e. so that for each time $t>0$ the function $p_\infty(\cdot,t)$ is the unique solution of the obstacle problem \eqref{eq:obstacle}.
We would now like to characterize the evolution of the saturated region. We start with the following proposition:
\begin{proposition}\label{prop:mu}
For all $t>0$, 
%let  $p(\cdot,t)$ be the unique solution of \eqref{eq:obstacle} (so that $p_\infty = p$ a.e.) 
%and 
$\mathcal P(t) :=\{p_\infty(\cdot,t)>0\}$ the positivity set of the solution of the obstacle problem \eqref{eq:obstacle}. 
%(recall that $p^+$ is defined as the trace of $p$, see \eqref{eq:tracep}).
Then the density equation in \eqref{eq:2} can be rewritten as
\begin{equation}\label{eq:rhomu} 
\pa_t \rho_\infty = \mu_t +\lambda\rho_{\infty}(1-\chi_{\mathcal{P}}),
\end{equation}
here $\mu_t:=\Delta p_\infty(\cdot,t) +\lambda(\cdot,t) \chi_{\mathcal{P}(t)}$, which is a non-negative Radon measure supported in  $\pa\mathcal{P}(t)\setminus \mathcal O(t)$.

\end{proposition}

When $\lambda\leq 0$, Equation \eqref{eq:rhomu} shows that the growth of $\rho_\infty$ can only occur when the measure $\mu$ is non zero (thus only on $\pa\mathcal{P}(t)\setminus \mathcal O(t)$) while  the density can only decay when $\rho_\infty(1 -\chi_{\mathcal{P}(t)})>0$. 
 Growth and decay thus take place according to different mechanisms. One is dictated by a singular measure, the other by an $L^\infty$ function. 
Note that $\mathcal{P}(t)$ is almost the saturated set $\Sigma(t)$, in the sense that their parabolic closures coincide (see Theorem~\ref{thm:velocity}.)

\medskip

Heuristically, \eqref{eq:rhomu} can have a geometric interpretation as follows.
Since $\rho_\infty=1$ in $\mathcal{P}(t)$, we can always write
$$ \rho_\infty(x,t) = \chi_{\mathcal{P}(t)}(x) + \rho^E(x,t) (1-\chi_{\mathcal{P}(t)}(x) ) $$
for some function $\rho^E$. Splitting the singular and regular part of \eqref{eq:rhomu}, we get the 
following:
\begin{equation}\label{eq:twoeq}
\begin{cases}
 (1-\chi_{\mathcal{P}(t)}(x) ) (\pa_t \rho^E - \lambda \rho^E) =0;\\
 (1-\rho^E(x,t)) \pa_t  \chi_{\mathcal{P}(t)} =\mu.
\end{cases}
\end{equation}
The first equation determines the value of $\rho_\infty$ outside of the congested set ($ \pa_t \rho^E = \lambda \rho^E$ when $p(x,t)=0$, supplemented by the condition that $\rho^E=1$ when $p(x,t)>0$).

\medskip

Formally, we have $\mu = |\nabla p|dS$, where $S$ is the surface measure on $\partial\mathcal{P}(t)$, so if $|\na p(x_0,t_0)|\neq0$, the second equation in \eqref{eq:twoeq} gives  $(1-\rho^E) V(x_0,t_0)=  |\na p(x_0,t_0)|$ (expansion of the congested region), while  if $|\na p(x_0,t_0)|=0 $, then either $ \pa_t  \chi_{\mathcal{P}(t)}=0$ or $\rho^E (x_0,t_0)=1$. The later  can only happen if 
$ \chi_{\mathcal{P}(t)}(x_0)=1$ as $t\to t_0^-$ and so 
$ \pa_t  \chi_{\mathcal{P}(t)}\leq 0$ (retraction of  the congested region). Altogether, this gives the free boundary condition \eqref{eq:fb}, assuming that the boundary of $\mathcal{P}(t)$ coincides with $\Sigma(t)$.

\medskip

Making these statement rigorous in the classical framework would require the development of a regularity theory which is not the topic of the present paper. 
Instead, in what follows, we will use the comparison principle to make sense of this in the spirit of viscosity solutions.

\medskip

\subsection{Velocity law: barrier approach}

%We begin with the comparison principle for weak solutions of the limit problem \eqref{eq:2}. 

%\begin{proposition}\label{uniqueness}
%Suppose $\lambda \in   L^2(0,T;H^1(\Omega))$. Let $(\rho,p)\in L^\infty(Q_T) \times L^2(0,T;H^1(\Omega))$ be a weak solution of \eqref{eq:2} in $\Omega \times [0,T]$. Then the following holds:
%\begin{itemize}
% \item[(a)] $(\rho,p)$ is unique such weak solution, and in particular is the unique limit of $(\rho_m, p_m)$ in $\Omega\times [0,T]$. \item[(b)] 
%For a domain $D\subset \Omega$ with smooth boundary, let
%$(\rho_1, p_1) \in L^\infty(D\times [t_1,t_2]) \times L^2(t_1,t_2;H^1(D))$ be a weak solution of 
% \begin{equation}\label{eq:eq3}
% \pa_t \rho= \Delta p - \lambda \rho \quad  \mbox{ in } \mathcal D'(D\times [t_1,t_2]), \qquad p \in P_\infty(\rho).
% \end{equation}
% If $p \leq  p_1$ on $\partial D \times [t_1,t_2]$ and $\rho \leq \rho_1$ on $t=t_1$, then  $\rho \leq \rho_1$ in $D \times [t_1, t_2]$. A corresponding result holds with the reverse inequalities.
% \end{itemize}
% \end{proposition}

We define
the external density $\rho^E$ in the set $\{p=0\}$ by solving the first-order equation $\pa_t \rho = \lambda \rho$ with appropriate boundary condition.
More precisely, given $x$, the open set $\mathrm{Int}( \{t \, ;\, p(x,t)=0\})$ can be written as $\cup_{i\in I} (a_i,b_i)$ and $\rho^E(x,t)$ for $t\in (a_i,b_i)$ is the solution of the first order ODE $\pa_t \rho = \lambda \rho$  with initial condition
$$ 
\rho^E(x,a_i) = 
\begin{cases}
 \rho_0(x) & \mbox { if } a_i=0 \\ 
1 & \mbox{ if } a_i>0
\end{cases}
$$

See Figure 1 above for an illustration.

\begin{figure}
\begin{tikzpicture}
\draw[red,ultra thick] (1,0) to  [out= 110, in=270] (0,1);
\draw[thick](0,1) to [out= 90, in= 270] (1.5,2);
\draw[red, ultra thick] (1.5,2) to [out=90, in = 270] (-1,3.5);
\draw[thick] (-1,3.5)  to [out=90, in =270] (2.5,5);
\draw[<->]  (-3,5.2) -- (-3,0) -- (6,0);
\node [left] at (-3,5.2) {$t$};
\node [below right] at (6,0) {$x$};
\draw[thick] (4,5) to [out=350, in=120] (4.5,3) to [out=50, in=270] (6,5);
\path [fill=lightgray] (1,0) to  [out= 110, in=270] (0,1) to [out= 90, in= 270] (1.5,2)
to [out=90, in = 270] (-1,3.5)  to [out=90, in =270] (2.5,5) --(-3,5)--(-3,0)--(1,0);
\path [fill=lightgray] (4,5) to [out=350, in=120] (4.5,3) to [out=50, in=270] (6,5) --(5,5);
\node at (-1,1) {$p>0$};
\node at (5, 4) {$p>0$};
\node at (3,3) {$p=0$};
\node  at (3,2) {$(\rho^E)_t=\lambda \rho^E$};
\draw[blue, ultra thick] (1,0)--(6,0);
\end{tikzpicture}
\caption{The external density $\rho^E$ has boundary values on red and blue parts of the boundary: $\rho^E=1$ on the red parts, and $\rho^E=\rho^0$ on the blue part.} 
\end{figure}
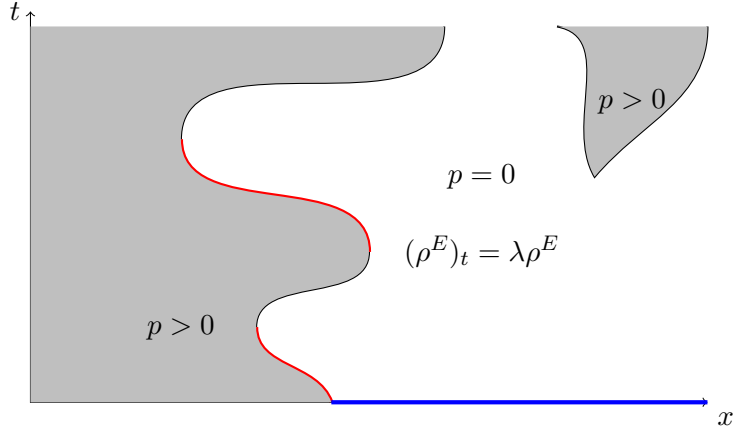

With this definition of $\rho^E$, and using the comparison principle for the limiting problem (Proposition \ref{prop:uniqueness}), we obtain the following description on the motion of the congested zone $\{\rho_{\infty}=1\}$:
\begin{theorem}\label{thm:velocity}
Let
$(\rho,p)\in L^\infty(Q) \times L^2_{loc}(0,\infty;H^1(\Omega))$ be a weak solution of \eqref{eq:2} with initial data $0\leq \rho^0(x)\leq 1$.  Then the following holds:
\begin{itemize}
\item[(a)] If $\lambda \in C(Q_T) \cap L^2(0,T;H^1(\Omega))$, then we have, in the sense of comparison with barriers, 
\begin{equation}\label{vel}
 (1-\rho^E)V= |\na p| \quad\hbox{ on } \partial\{\rho=1\}.
 \end{equation}
%Here the external density $\rho^E: \overline{\{p=0\}} \to [0,1]$ solves $\rho_t =  \lambda \rho$ in $(\tau(x,t),t)$  \textcolor{black}{where $\tau(x,t)$ was the most recent time before $t$ that $x$ was in the support of $p$, i.e. $\tau(x,t) := \inf\{0\leq\tau \leq t : (x,s)\in \overline{\{p=0\}} \mbox{ for } \tau \leq s\leq t\}$,  with initial data $1$ if $\tau(x,t)>0$, and otherwise with initial data $\rho^0(x)$.  }
\item[(b)] $\rho^E$ coincides with $\rho$ a.e. outside of $\overline{\{\rho=1\}}$.
\item[(c)] If $\lambda$ is negative, then for any $T>0$
$$
\overline{\{p>0\} \cap Q_T}= \overline{\{p>0\}}\cap Q_T= \overline{\{\rho=1\}}\cap Q_T.
$$
\end{itemize}
\end{theorem}

The barriers used to make sense of \eqref{vel} are local smooth sub- and super-solutions of \eqref{eq:2}. Their description can be found in in section  \ref{sec:velocity}.  Such comparison property is akin to the viscosity solutions approach taken by \cite{KP} for $\lambda>0$. We do not touch upon the issue of whether the barrier properties are enough for a complete characterization of the limit solution: see \cite{BKP} and \cite{KP} for analysis in this direction.

\medskip

Part (c) in above theorem says that when $\lambda$ is negative, the closure of the pressure support coincides with that of $\{\rho=1\}$, and that the congested zone $\{\rho=1\}$ cannot all of a sudden expand.  This is not true when $\lambda$ is positive, due to the nucleation of the congested zone generated by the growth of the external density. The set $\{\rho=1\}$ certainly can discontinuously shrink. For instance if $\lambda$ decreases over time, the pressure decreases and the set $\{\rho=1\}$ may start shrinking. While shrinking, if a component of the set gets disconnected at $t=t_0$ from $K$, the pressure in this region will drop to zero and $\rho$ will immediately decrease below one after $t_0$. Such scenario makes it difficult to describe $\rho^E$ in an explicit way, except when $\lambda$ only increases over time.

\begin{theorem}
Suppose that $\lambda \in C(\Omega \times [0,T]) \cap L^2(0,T;H^1(\Omega))$ is non-decreasing over time, and let $(\rho,p)$ be the weak solution of \eqref{eq:2} with initial data $\rho^0\in BV$.  Then the set
 $\{p(\cdot,t)>0\}$ is monotone increasing in time. Moreover for all $t\geq 0$
$$
 \rho(\cdot,t) = \chi_{\Sigma(t)} + \rho^E\chi_{\RR^n\setminus \Sigma(t)}, \hbox{ where } \rho^E(x,t):= \rho^0(x)\exp^{\int_0^t \lambda(x,s) ds}.
$$
In particular $\overline{\Sigma(t)} = \overline{\{\rho(\cdot,t)=1\}}$ for all $t>0$. 

\medskip

If $\rho^0$ is a characteristic function and $\Sigma_0 = \{\rho^0=1\} = \{\rho^0>0\}$, then $\rho$ remains a characteristic function for all positive times.
\end{theorem}
Note that we may initially have $\{p(\cdot,0)>0\}$ as a strict subset of  $\{\rho^0=1\}.$ In this case this last theorem states that $\{\rho=1\}$ experiences an initial discontinuous shrinkage.

\subsection{Numerical examples}

\begin{figure}\label{fig1}
\begin{tabular}{ccc}
\scalebox{.18}{\pdfimage{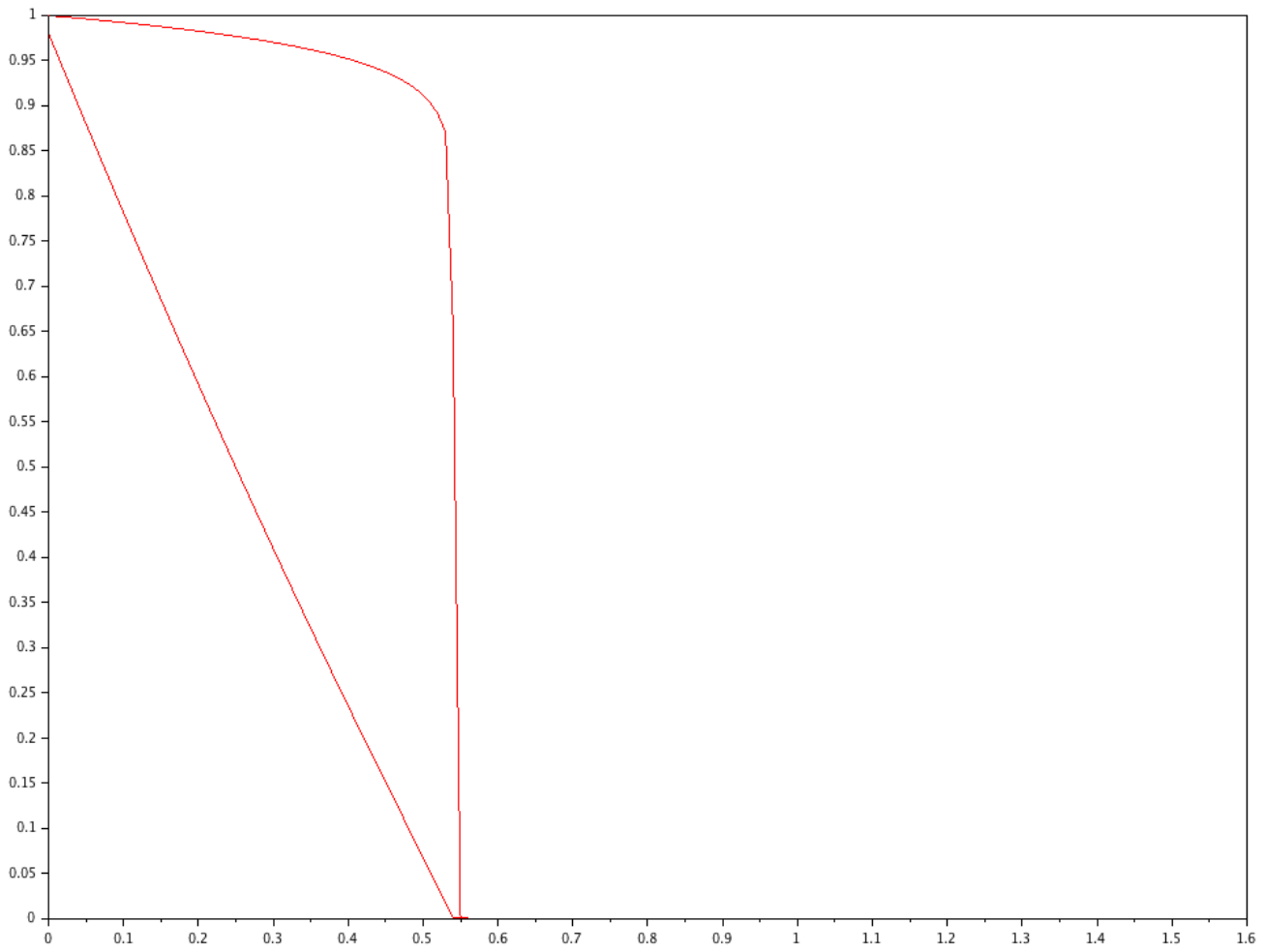}} & \scalebox{.18}{\pdfimage{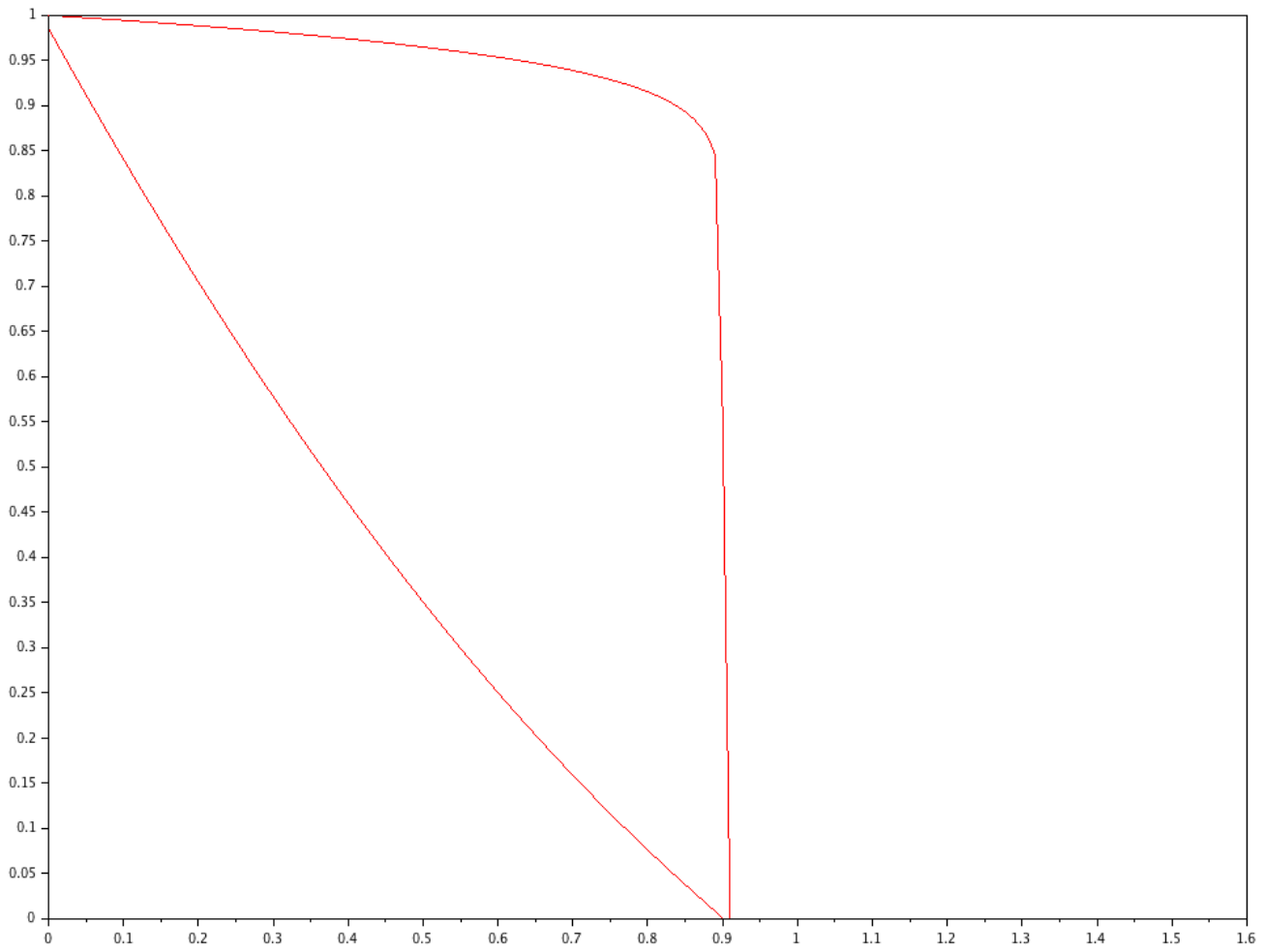}} & \scalebox{.18}{\pdfimage{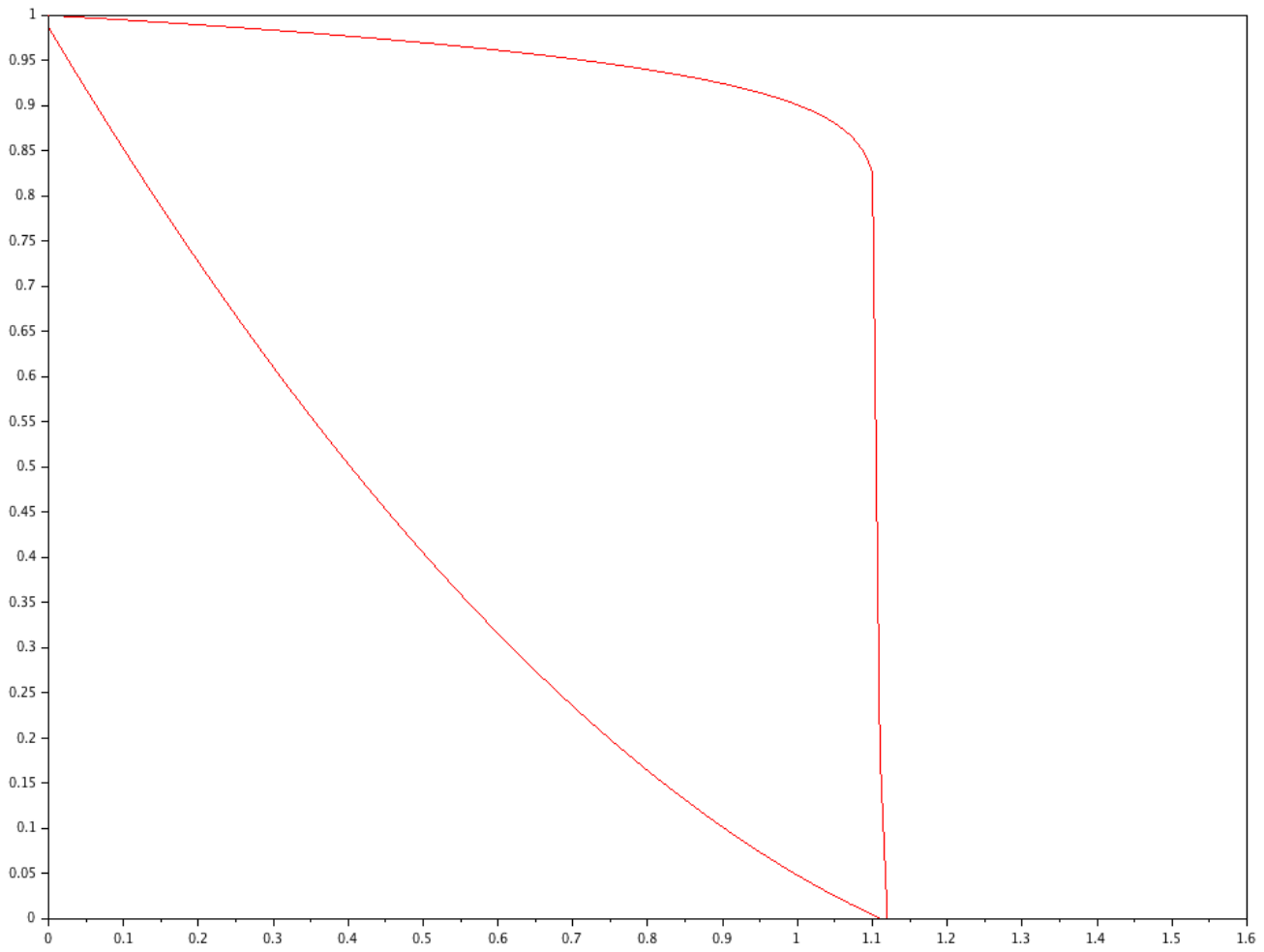}} \\ 
\scalebox{.18}{\pdfimage{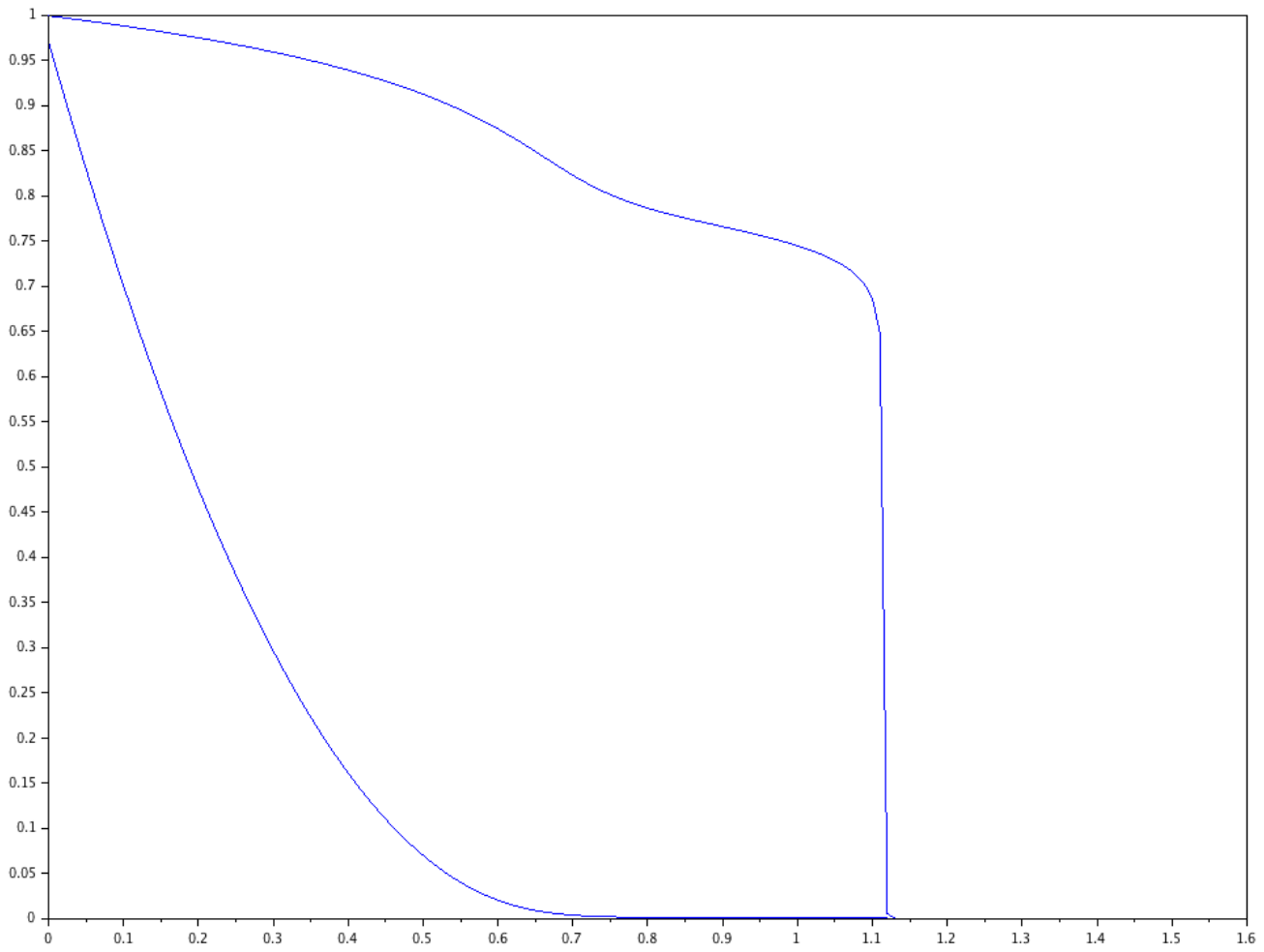}}& 
\scalebox{.18}{\pdfimage{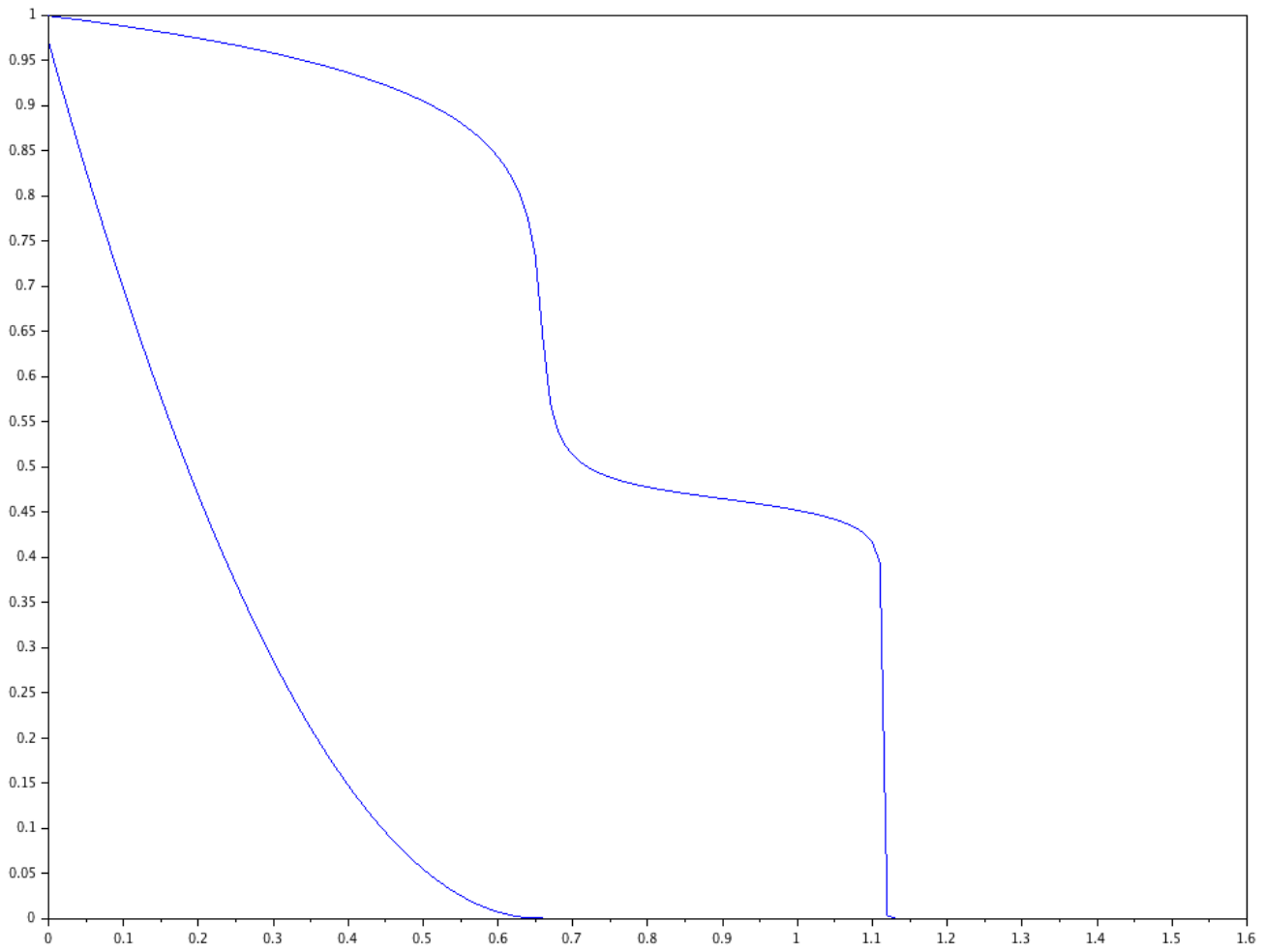}} & \scalebox{.18}{\pdfimage{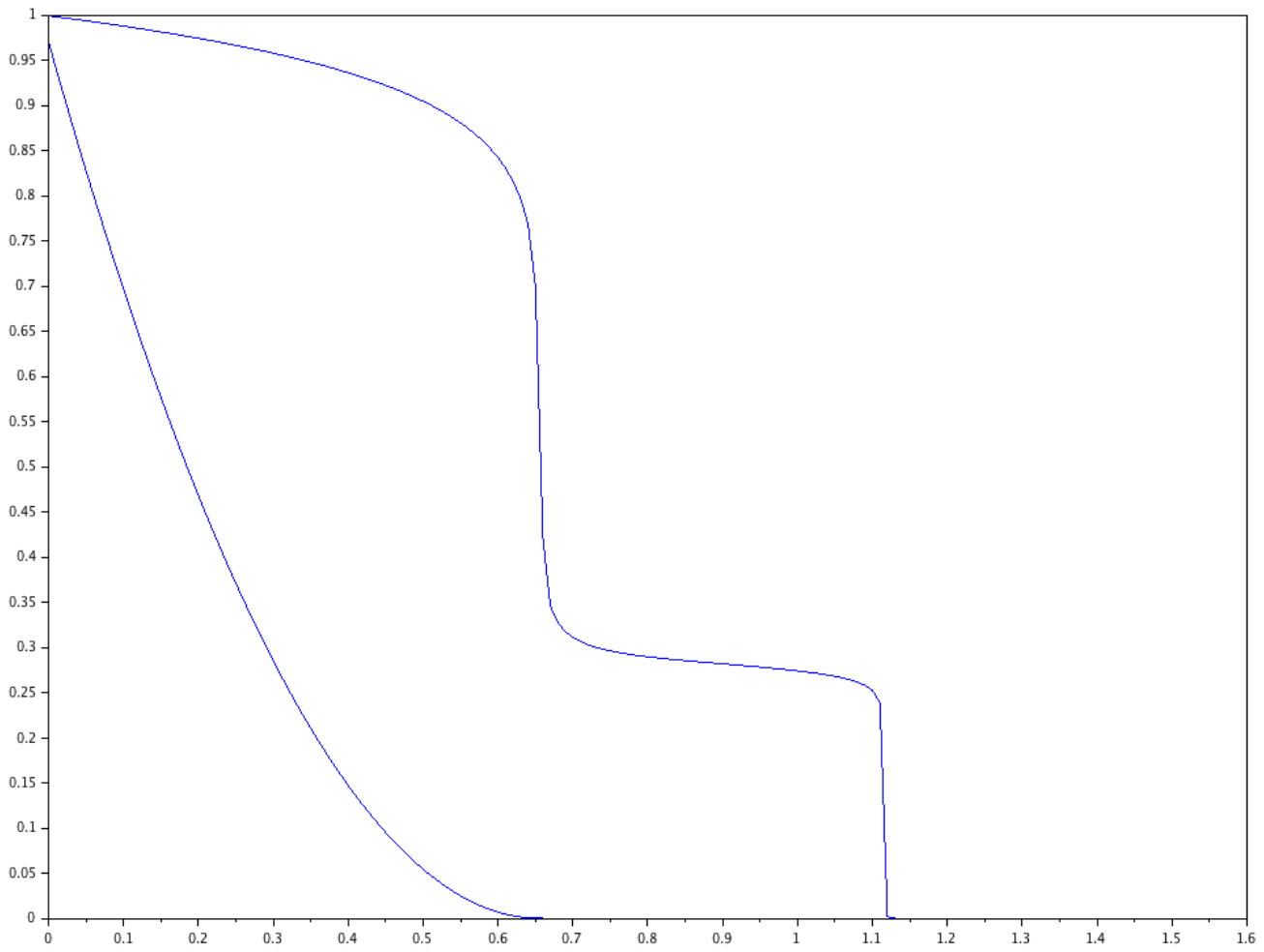}} \\
 \scalebox{.18}{\pdfimage{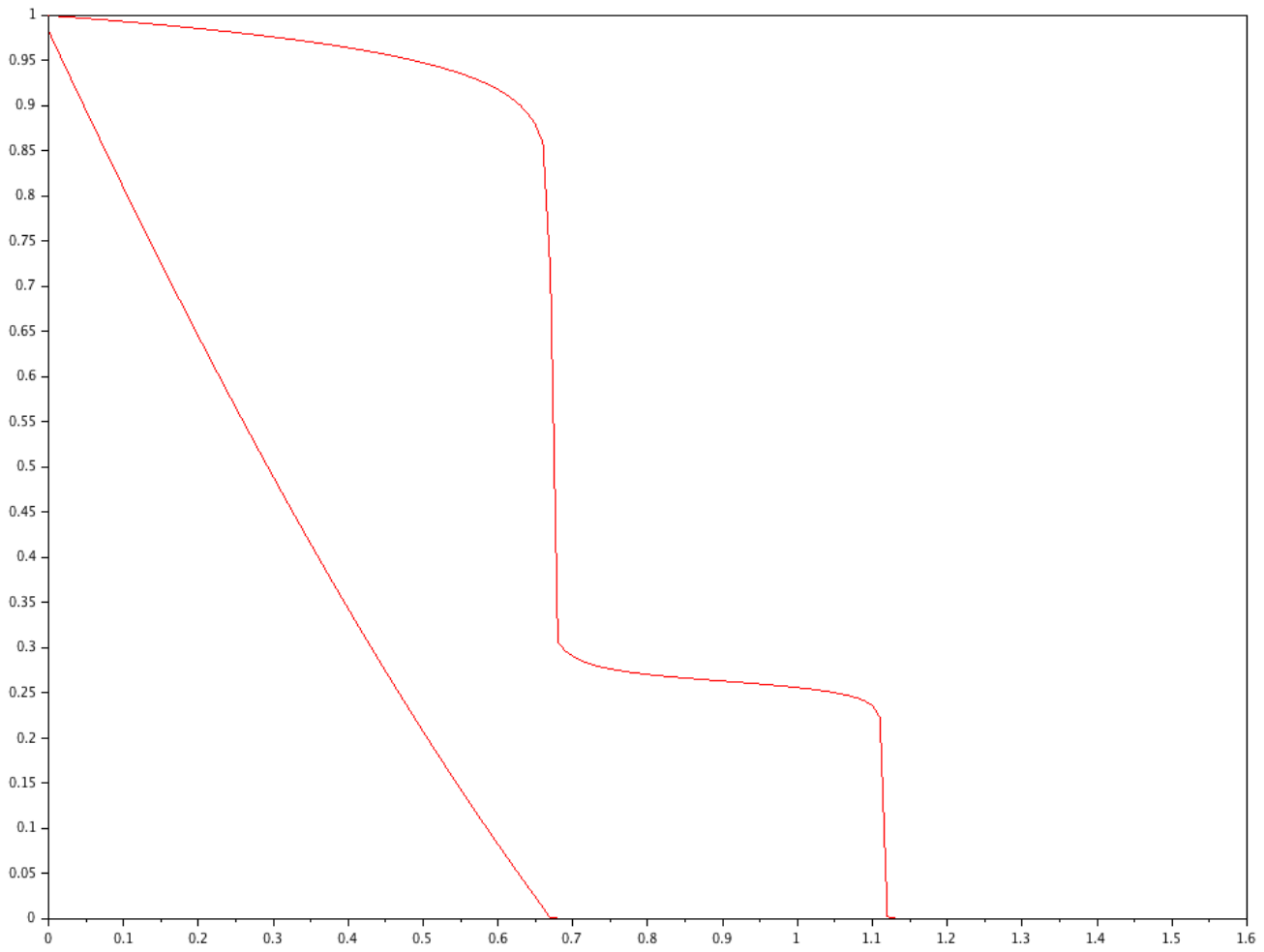}} & \scalebox{.18}{\pdfimage{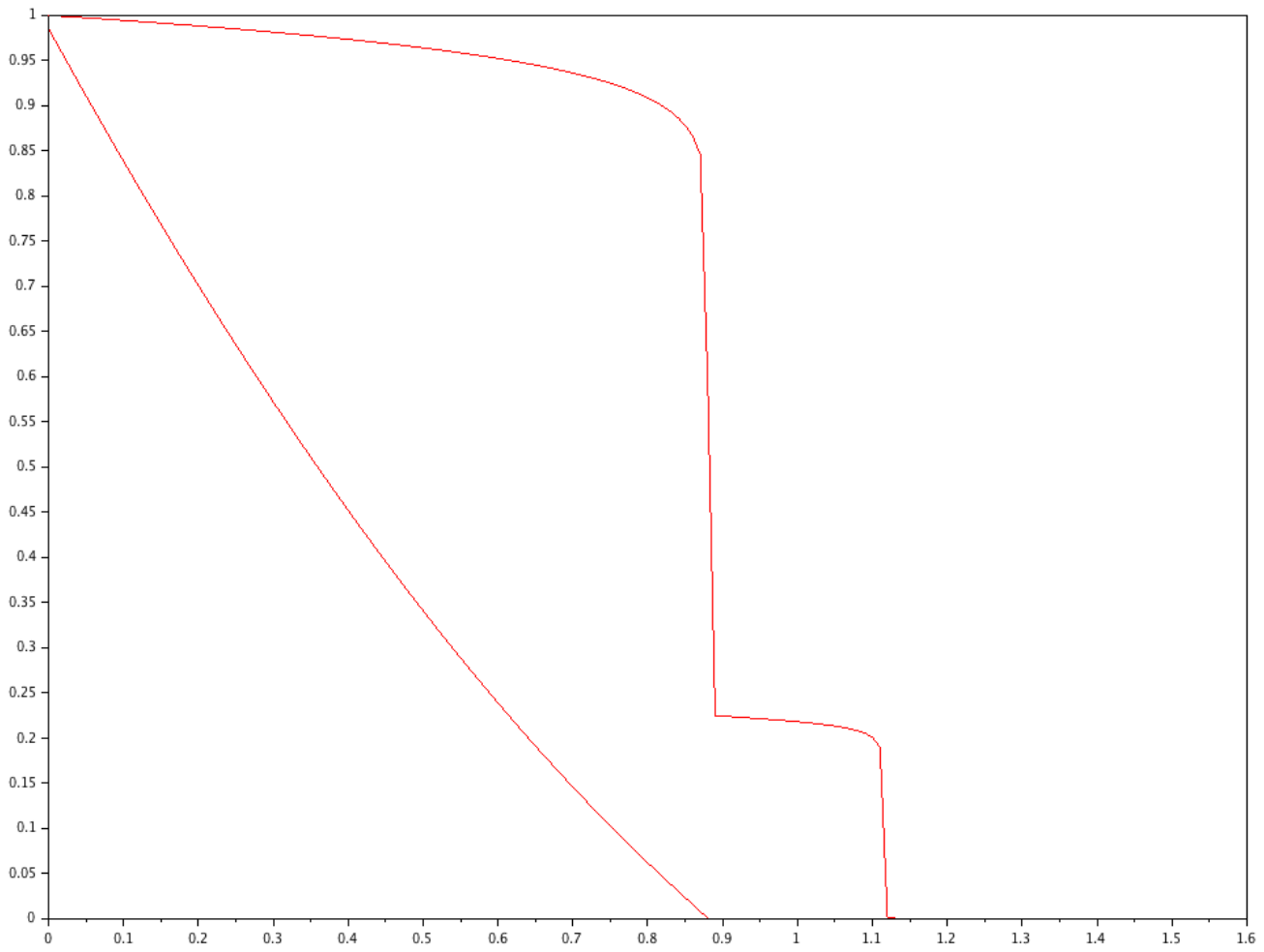}}  &\scalebox{.18}{\pdfimage{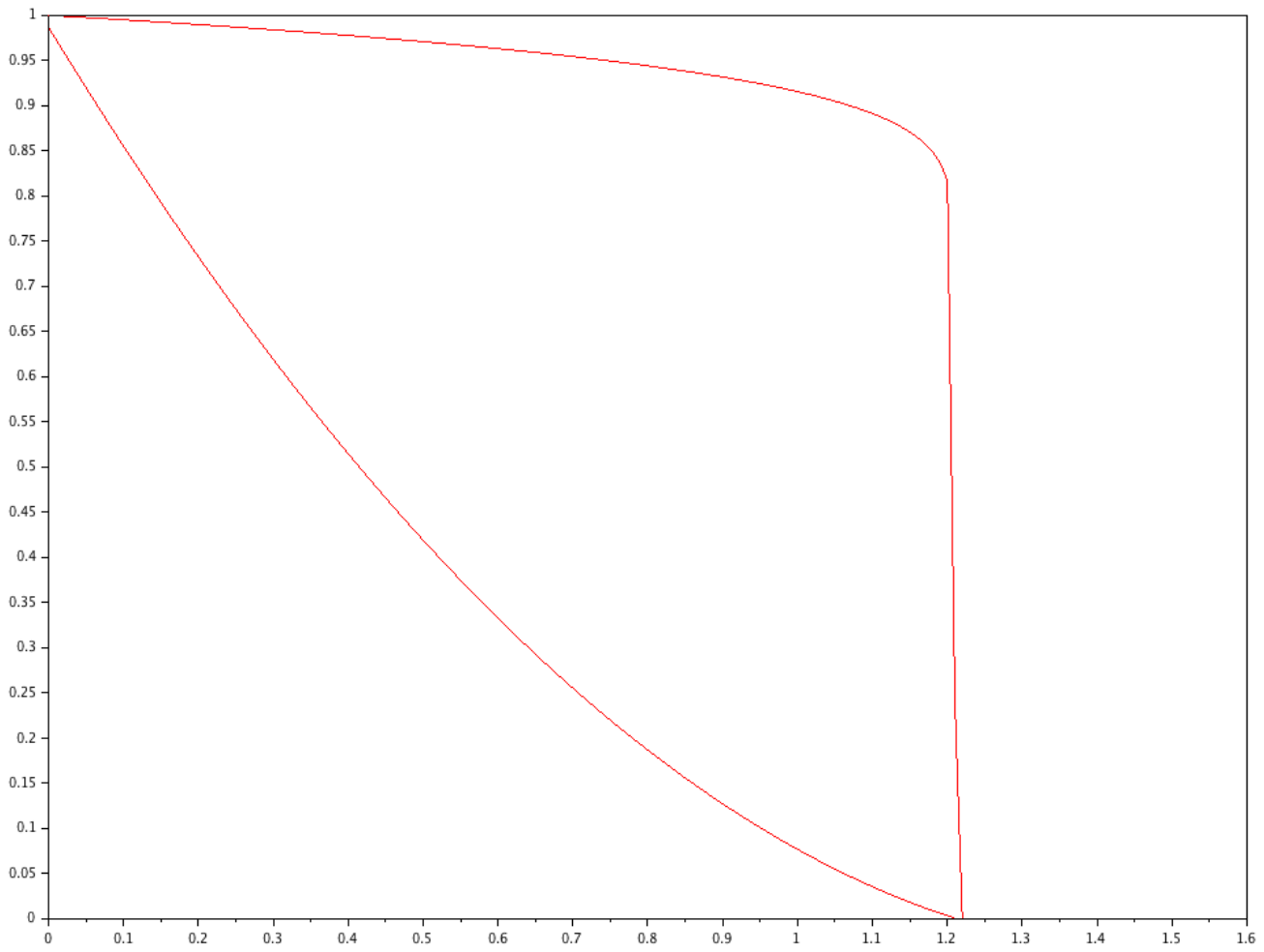}}
\end{tabular}
\caption{Graph of the density (upper curve) and pressure (lower curve) when $m=40$ and 
$\lambda$ given by \eqref{eq:lambdafig}.
The functions are shown at time $t=0^+$, $t=0.35$, $t=0.75^-$, $t=0.75^+$, $t=0.9$, $t=1^-$, $t=1^+$, $t=1.2$ and $t=1.8$}
\end{figure}

Figure \ref{fig1} shows the evolution of the density and pressure in a simple framework to illustrate the receding and expanding motion of the free boundary.
We consider the one dimension porous media equation
$$
 \pa_t \rho- \pa_x(\rho\pa_x p) = \lambda(t) \rho \qquad  \mbox{ in }  (0,\infty)\times (0,T), \qquad p = \frac{m}{m-1} \rho^{m-1}
$$
with the boundary condition $\rho(0) =1$ and $m=40$ (so we are close to the limiting problem. In particular, the density is close to, but not equal to $1$ when $p>0$).
The coefficient $\lambda(t)$ is independent of $x$ but changes value discontinuously in time:
\begin{equation}\label{eq:lambdafig}
\lambda (t) = 
\begin{cases}
 -1& \mbox{ if } t\in [0,.75)\\
  -5 & \mbox{ if } t\in [.75,1)\\
  -1 &  \mbox{ if } t\geq 1.
\end{cases}
\end{equation}
The set $\{p(t)>0\}$ is expanding with finite speed for $t\in (0,.75)$ (first row) and receding instantaneously at $t=.75^+$. The density is then decreasing for $t\in [.75,1)$ in the region where $p=0$ since $\pa_t \rho = -5\rho$ in that region (second row).
Finally, for $t>1$ (third row) the  set $\{p(t)>0\}$ is again expanding with finite speed.

\medskip

\section{Proof of Theorem \ref{thm:1}}
The proof of this theorem uses many classical techniques (see in particular \cite{PQV}), though we have to be careful with the two main differences between our framework and that of \cite{PQV}: the lack of  sign of $\lambda$ and the presence of the fixed boundary $\pa K$.

\subsection{Notion of solutions for \eqref{eq:1}}

First, we recall some well known facts about the porous media equation \eqref{eq:1} (we refer the interested reader to \cite{Vaz2007}, Chapters 5 (Definition 5.5 and Theorem 5.14).
\begin{definition}
For  $\rho^0\in L^1(\Omega)$, $g \in L^2(0,T;H^1(\Omega))$ and $\lambda\in L^1(Q_T)$, we say that a non-negative function $\rho\in L^1(Q_T)$ is a {\it weak solution} of \eqref{eq:1} with $\rho^0_m=\rho^0$ and $f:= g^{1-1/m}$ if
\begin{itemize} 
\item[(i)] $\rho^m \in L^2(0,T;H^1(\Omega))$ with its trace on $\partial K \times [0,T]$ equal to $g$;
 \item[(ii)] $\rho\in L^2(Q_T)$;
 \item[(iii)] $\rho$ satisfies the identity
 $$
\int\int_{Q_T} (\rho\pa_t \psi -\nabla \rho^m \cdot \nabla \psi +\lambda \psi ) dx dt = - \int_{\Omega} \rho_m^0(x)\psi(x,0) dx
$$
for any function $\psi\in C^1(\overline{Q}_T)$ which vanishes on $\partial K \times [0,T]$ and for $t=T$. 
\end{itemize}
\end{definition}    
Existence of a weak solution can be established by approximation with smooth functions, which  either solves the porous media equation with strictly positive initial data or solves a regularized equation with strictly positive diffusion (see Theorem 5.14 of \cite{Vaz2007}). Uniqueness of the weak solution is a consequence of the following comparison principle, which we will use often in our analysis. 
 \begin{lemma}\label{cp:pme} 
Let $ \rho$ and $\tilde{\rho}$ be two weak solutions of \eqref{eq:1} with initial data $\rho^0_m, \tilde{\rho}^0_m$ and fixed boundary data $f$ and $\tilde{f}$.  If $\rho^0_m \leq \tilde{\rho}^0$ a.e. and $f \leq \tilde{f}$ a.e., then $\rho \leq \tilde{\rho}$ a.e. 
\end{lemma}

\subsection{Maximum principle: $L^\infty$ bounds for $\rho_m$ and $p_m$ and $\na p_m\cdot \nu|_{\pa K}$}

\begin{lemma}\label{lem:basic}
Under conditions 
\eqref{eq:lambdabd}, \eqref{eq:condbc}
and
\eqref{eq:condin1}, and for all $T>0$,
there exists a constant $C=C(T)>0$ independent of $m$ such that the following holds: 
\item For sufficiently large $m$ (depending on $T$) the pressure $p_m$ satisfies:
\begin{equation}\label{eq:infty}
 0\leq p_m(x,t)\leq C \qquad  \mbox{ for all } (x,t)\in Q_T ,
\end{equation}
and
\begin{equation}\label{eq:flux}
-C \leq \nu(x)\cdot \na p_m(x,t)\leq C \quad \mbox{ for all } x\in \pa K, \quad 0 \leq t \leq T.
\end{equation}
\item Moreover,
\begin{equation}\label{eq:rho1}
\rho_m(x,t) \to 1,\quad\mbox{locally uniformly in $U\times \RR_+$}
\end{equation}
for some neighborhood $U$ of $K$ and
\begin{equation}\label{eq:supprho}
 \supp \rho_m(\cdot,t) \subset B_{\overline R+C(T)} \qquad \mbox{ for all $t\in(0,T)$.}
 \end{equation}
\end{lemma}
\begin{remark}
Note that by \eqref{eq:rho1}, $\rho_m$ stays uniformly positive and solves a uniformly parabolic equation in $U$. 
It is thus smooth, a  fact we will use repeatedly when dealing with the boundary data on $\pa K$.
\end{remark}
\begin{proof}
We fix $T>0$.
This lemma follows from the maximum principle for the pressure $p_m$, which, we recall solves
$$\pa_t p = (m-1)p (\Delta p +\lambda) + |\na p|^2.$$
In view of \eqref{eq:uphi}, $\underline \vphi(x)$ satisfies $\Delta \underline \vphi + \lambda = \Lambda+\lambda \geq 0$ and is therefore a subsolution for this equation. 
Assumption \eqref{eq:condin1} thus implies
$$ p_m(x,t) \geq \underline \vphi(x) \qquad \forall (x,t)\in\Omega\times\RR^+.$$

For the upper barrier, we define the function $\bar{u}(x,t)$ as follows:
For all $t>0$, the function 
$x\mapsto  \bar{u}(x,t)$ solves 
$$
-\Delta v = \Lambda +1 \mbox{ in } B_{R(t)} \setminus K, \quad v = f^\frac{m}{m-1} \mbox{ on } \pa K,\quad v= 0 \mbox{ on } \pa B_{ R(t)}
$$
where $R(t):= \ds \overline R+ \int_0^t M(s) ds$, with $M(s):= 2\sup_{x\in \pa B_{R(t)}} |\nabla\bar{u}(\cdot,t)|$. The function $\bar u$  is extended by $0$ outside $ B_{ R(t)}$.
Since $R(t)$ depends on $\bar u(x,t)$, the function $\bar{u}$ can be constructed for instance by discrete-time approximation. 
We note that \eqref{eq:ophi} implies in particular that $\bar u(x,0) =\overline \vphi(x)$

We  claim then $\bar{u}$ is a supersolution for the pressure equation  for sufficiently large $m$. To see this, note first that when $\bar u \geq (m-1)^{-1/2}$ we have
$$
\pa_t \bar u  \geq 0 \geq (m-1)\bar u(\Delta \bar u + \lambda) + |\nabla \bar u|^2 \quad\,\, \hbox{ if }  m \geq \sup_{0\leq t\leq T} |\nabla  \bar u|^4(\cdot,t).
$$
On the other hand,  since $\pa_t \bar u \geq 2|\na \bar u|^2>0$ on its zero level set $\pa B_{R(t)}$, it is clear that for small enough  $\epsilon=\epsilon(T)$  we have $\pa_t \bar u \geq |\na \bar u|^2$ in $0\leq \bar u \leq \e$. Our  claim follows if $m$ is large enough (depending on $T$).

%, and, by comparison principle for the pressure equation (which follows from Lemma~\ref{cp:pme}), \eqref{eq:infty} follows.

%\medskip

%For the lower barrier, we can just  take the stationary profile $\underline \vphi$ from \eqref{eq:uphi} for the pressure equation. 
The comparison principle for the pressure equation now yields:
$$
\underline \vphi(x) \leq p_m(x,t)\leq \bar{u}(x,t) \qquad\forall (x,t)\in Q_T.
$$

The results now follow:
\eqref{eq:infty} follows from upper bound, while
the lower bound together with the fact that  $\rho_m \sim p_m^{\frac 1 {m-1}}$ implies \eqref{eq:rho1}.  
The fact that $\bar u$ is supported in $B_{R(t)}$  implies \eqref{eq:supprho}
and since $\underline \vphi(x)=  p_m(x,t) = \bar{u}(x,t)$ on $\pa K$, we get
$$
-C \leq \nu(x)\cdot \na \underline \vphi(x) \leq \nu(x)\cdot \na p_m(t,x) \leq \nu(x)\cdot \na\bar{u}(x,t)  \leq C \qquad \forall x\in \pa K, \; t\in (0,T].
$$

 \end{proof}

\subsection{$L^1$ bounds for $\rho_m$ and $p_m$}
\begin{lemma}\label{lem:L1}
For all $T>0$, 
there exists a constant $C(T)$ depending on $\Lambda$ and $T$ such that
\begin{equation}\label{eq:L1rho}
\| \rho_m(t) \|_{L^1(\Omega)} \leq \| \rho^0_m \|_{L^1(\Omega)}  e^{\Lambda T}  +C (T) 
\end{equation}
and
\begin{equation}\label{eq:L1p}
\| p_m(t) \|_{L^1(\Omega)} \leq  C \| \rho^0_m \|_{L^1(\Omega)}  e^{\Lambda T}  +C (T)
\end{equation}
for $t\in [0,T]$ and  $m\geq 2$.
\end{lemma}

\begin{proof}
Integrating \eqref{eq:1} on $\Omega$ yields
\begin{align*}
\frac{d}{dt}\int_{\Omega} \rho_m(t)\, dx
& =  \int_{\Omega} \lambda(t)\rho_m(t) \,dx + \int_{\pa K} \rho_m \na p_m\cdot\nu\, dS\\
& \leq  \Lambda  \int_{\Omega} \rho_m(t) \,dx + C,
\end{align*}
where we used  \eqref{eq:infty}, \eqref{eq:flux}. 
The bound \eqref{eq:L1rho} follows by a Gronwall argument.
The bound \eqref{eq:L1p} then follows from \eqref{eq:L1rho} and \eqref{eq:infty} since  $p_m= \frac{m}{m-1}\rho_m^{m-1}\leq \frac{m}{m-1}C \rho_m$.
\end{proof}

\subsection{Bounds on the derivatives of $\rho_m$ and $p_m$}
For $\delta>0$, we define
$$\Omega^\delta := \{x\in \RR^n\,;\, \mathrm{dist}(x,K)>\delta\} .$$
\begin{lemma}\label{lem:BV}
%Assume that 
%\begin{equation}\label{eq:condinitt}
%\| \Delta p_{in} -\lambda \rho^0_m\|_{L^1(\Omega)} + \|\na  \rho^0_m   \|_{L^1(\Omega)} \leq C
%\end{equation}
%for some constant $C$ independent on $m$. 
For any $\delta>0$, there exists a constant $C_\delta$ independent on $m$ such that 
\begin{align}
\|\pa_t \rho_m (t) \|_{L^1(\Omega^\delta)} 
 \leq C_\delta    & \qquad \forall t>0\label{eq:BVtrho} \\
\|\pa_{x_i} \rho_m (t) \|_{L^1(\Omega^\delta)} 
\leq C_\delta    & \qquad \forall t>0 \label{eq:BVxrho} 
\end{align}
Similarly, denoting by $B_R$ the ball of radius $R$, we have the following bounds:
\begin{align}
\|\pa_t p_m \|_{L^1((0,T)\times \Omega^\delta\cap B_R)} 
& \leq C_{\delta,R,T}  \label{eq:BVtp} \\
\|\pa_{x_i} p_m \|_{L^1((0,T)\times \Omega^\delta\cap B_R)} 
& \leq C_{\delta,R,T}   \label{eq:BVxp} 
\end{align}

\end{lemma}

\begin{proof}
Proceeding as in \cite{PQV}, we differentiate the first equation in \eqref{eq:1} 
with respect to time and multiply it by $\mathrm{ sign }(\pa_t \rho_m)$ and use Kato's inequality to obtain 
\begin{equation}\label{eq:patrho}
\pa_t |\pa_t \rho_m| - \Delta( m \rho_m^{m-1} |\pa_t \rho_m| ) \leq \lambda |\pa_t \rho_m| +\rho_m| \pa_t \lambda| 
\quad \mbox{ in } \Omega.
\end{equation}
We cannot simply integrate this equation over $ \Omega$ because of the boundary condition on $\pa K$. Instead, given a large ball $B_R$ such that $K\subset B_R$, we introduce the function $\vphi$ such that
$ \vphi = 0$ on $\pa K$, $\vphi=1$ on $\pa B_R$, $\Delta \vphi = 0$ in $B_R\setminus K$ and we extend this function by $1$ outside $B_R$. This function satisfies
$$ \vphi|_{\pa K} = 0, \quad \Delta \vphi \leq 0 \mbox{ in } \Omega, \quad \vphi>0 \mbox{ in } \Omega.$$
Multiplying \eqref{eq:patrho} by $\vphi$ and integrating over $\Omega$, and using the fact that $\vphi|_{\pa K} = 0$ and $m \rho_m^{m-1} |\pa_t \rho_m|   = \frac{m}{m-1} f^\frac{1}{m-1} \pa_t f $ on $\pa K$, we deduce
\begin{align*}
   \frac {d}{dt} \int_{\Omega} |\pa_t \rho_m| \vphi\, dx 
  & \leq  \int_{\Omega}   m \rho_m^{m-1} |\pa_t \rho_m| \Delta \vphi\, dx
  -\int_{\pa K} m \rho_m^{m-1} |\pa_t \rho_m| \na \vphi\cdot \nu\, dS \\
&\quad     +\int_{\Omega}  \lambda |\pa_t \rho_m| \vphi\, dx   + \int_{\Omega}   \rho_m | \pa_t \lambda | \vphi\, dx\\
    &\leq C+\Lambda  \int_{\Omega} |\pa_t \rho_m| \vphi\, dx
  +  \int_{\Omega}   \rho_m | \pa_t \lambda | \vphi\, dx.
  \end{align*}
Since $\pa_t \rho_m (0)= \Delta \rho^{m}_{in}  + \lambda  \rho^0_m $, the bound \eqref{eq:condin3}  implies $ \|\pa_t \rho_m (0)   \|_{L^1(\Omega)} \leq C$ and using \eqref{eq:supprho}, we deduce 
$$
\|\pa_t \rho_m (t) \vphi \|_{L^1(\Omega)} \leq C(T)  +  C \int_0^T e^{\Lambda(T-t)}\int_{B_{R_0+CT}} |\pa_t \lambda(x,s)|\, dx\, dt
$$
%$$
%\|\pa_t \rho_m (t) \vphi \|_{L^1(\Omega)} \leq C(T)
%+   \|\pa_t \rho_m (0)   \|_{L^1(\Omega)} +   \int_0^t  \int_{\Omega}   | \rho_m(x,t)| |\pa_t \lambda(x,t)| \, dx\, ds
%$$
%Since $\pa_t \rho_m (0)= \Delta \rho^{m}_{in}  - \lambda  \rho^0_m $, the bound \eqref{eq:condin3}  and \eqref{eq:supprho} imply
%$$
%\|\pa_t \rho_m (t) \vphi \|_{L^1(\Omega)} \leq C  +  C \int_0^T \int_{B_{R_0+CT}} |\pa_t \lambda(x,s)|\, dx\, dt
%$$
and \eqref{eq:BVtrho} follows from \eqref{eq:lambdaBV} and the fact that $\min_{\Omega^\delta} \vphi>0$ for all $\delta>0$ (by the strong maximum principle).
\medskip

To get an estimate on $\pa_t p_m$, we want to take advantage of the term $\int_{\Omega}   m \rho_m^{m-1} |\pa_t \rho_m| \Delta \vphi\, dx$ in the inequality above.
We thus define, for $\eta>0$,
$\vphi_\eta$ such that
$ \vphi_\eta = 0$ on $\pa K$, $\vphi_\eta=1$ on $\pa B_R$, $\Delta \vphi_\eta = -\eta $ in $B_R\setminus K$ and we extend this function by $1$ outside $B_R$. 

Given $R$,  we claim that $\Delta \vphi_{\eta} \leq 0$ in $\Omega$ if $\eta$ is sufficiently small (depending on $R$). 
Indeed, Hopf's  Lemma implies $x\cdot \nabla \vphi_0 >0$ on $\pa B_R$, so the $C^1$-convergence of $\vphi_{\eta}$ to $\vphi_0$ implies $x \cdot \nabla \vphi_{\eta} \geq 0$ on $\pa B_R$. 

Proceeding as above, we get:
\begin{align*}
   \frac {d}{dt} \int_{\Omega} |\pa_t \rho_m| \vphi_\eta\, dx 
   + 4\eta \int_{\Omega}   m \rho_m^{m-1} |\pa_t \rho_m| \, dx\leq C(T)
\end{align*}
Integrating in $t$, we deduce that for all $\delta>0$, $R>0$ and $T$, there exists $C_{\delta,R,T}$ such that
$$
\int_0^T \int_{\Omega^\delta\cap B_R}   m \rho_m^{m-1} |\pa_t \rho_m|  \, dx\, dt \leq C_{\delta,R,T}.$$
Finally, we write
\begin{align*}
\int_0^T \int_{\Omega^\delta\cap B_R}  |\pa_t p_m| \, dx\, dt  
& =
\int_0^T \int_{\Omega^\delta\cap B_R} m \rho_m^{m-2} |\pa_t \rho_m| \, dx\, dt \\
&  \leq \int_0^T \int_{\Omega^\delta\cap B_R\cap\{\rho_m<1/2\}} m (1/2)^{m-2} |\pa_t \rho_m| \, dx\, dt \\
& \qquad+ 2 \int_0^T \int_{\Omega^\delta\cap B_R\cap\{\rho_m>1/2\}} m \rho_m^{m-1} |\pa_t \rho_m| \, dx\, dt \\
&  \leq \int_0^T \int_{\Omega^\delta} |\pa_t \rho_m| \, dx\, dt + C_{\delta,R,T}
\end{align*}
which gives \eqref{eq:BVtp}.

\medskip

We proceed similarly for the bound on $\pa_{x_i}\rho_m$. The only difference is that we do not have $\pa_{x_i} \rho_m |_{\pa K} = 0$, so we have an additional boundary term to worry about.
More precisely, differentiating the first equation in \eqref{eq:1} 
with respect to $x_i$, multiplying it by $\mathrm{ sign }(\pa_{x_i} \rho_m)$ and using Kato's inequality, we obtain
\begin{equation}\label{eq:paxrho}
\pa_t |\pa_{x_i} \rho_m| - \Delta( m \rho_m^{m-1} |\pa_{x_i} \rho_m| ) \leq \lambda |\pa_{x_i} \rho_m| + |\pa_{x_i} \lambda| \rho_m
\quad \mbox{ in } \Omega.
\end{equation}
With the same cut-off function $\vphi$ as above, we get
\begin{align*}
\frac {d}{dt} \int_{\Omega} |\pa_{x_i} \rho_m| \vphi\, dx 
& \leq 
 \int_{\Omega}   m \rho_m^{m-1} |\pa_{x_i} \rho_m| \Delta \vphi\, dx 
 - \int_{\pa K} m \rho_m^{m-1} |\pa_{x_i} \rho_m|  \na \vphi\cdot \nu\, dS \\
&\qquad + \int_{\Omega}  \lambda |\pa_{x_i} \rho_m| \vphi\, dx+ \int_{\Omega}   | \rho_m | |\pa_{x_i} \lambda| \, dx  \\
& \leq   \int_{\pa K}  \rho_m  |\pa_{x_i} p_m|  |\na \vphi\cdot \nu|\, dS +\Lambda  \int_{\Omega}    |\pa_{x_i} \rho_m| \vphi\, dx +\int_{\Omega}   | \rho_m | |\pa_{x_i} \lambda| \, dx 
\end{align*}
To conclude, we thus note that the estimate \eqref{eq:flux} gives a bound on the normal derivative of $p_m$ on $\pa K$, while the condition $p_m|_{\pa K} = \frac{m}{m-1} f $ together with the regularity assumptions \eqref{eq:condbc} implies that the tangential derivatives of $p_m$ are uniformly bounded on $\pa K$. We deduce
that $ |\pa_{x_i} p_m| |_{\pa K}\leq C$ and so (using  \eqref{eq:condin3}, \eqref{eq:supprho} and \eqref{eq:lambdaBV}):
$$ \frac {d}{dt} \int_{\Omega} |\pa_{x_i} \rho_m| \vphi\, dx \leq C(T) \qquad \forall t\in(0,T). $$
Hence
$$
\|\pa_{x_i} \rho_m (t) \vphi \|_{L^1(\Omega)} \leq   \|\pa_{x_i} \rho_m (0)   \|_{L^1(\Omega)} + C(T),
$$
and \eqref{eq:BVxrho} now  follows from  \eqref{eq:condin3}.

\end{proof}

\subsection{Passing to the limit}
We denote
$$ \Omega_k := \{x\in \Omega\,;\, \mathrm{dist}(x,K) > 1/k, \; |x|\leq k\}.
$$
Lemma \ref{lem:BV} implies that $\rho_m$ and $p_m$ are bounded in $\mathrm{BV}(\RR_+\times \Omega_k)$ for all $k$ and thus converge (up to a subsequence) strongly  in $L^1([0,k]\times \Omega_k)$. By a diagonal extraction process, we can thus find subsequences (still denoted $\rho_m$ and $p_m$) and functions $\rho_\infty$, $p_\infty$ such that
$\rho_m$ (resp. $p_m$) converges to $\rho_\infty$ (resp. $p_\infty$) strongly in $L^1_{loc} (\RR_+\times \Omega)$.

\medskip

Next, we note that 
$$ \rho_m \, p_m = \left(\frac{m-1}{m}\right) ^\frac{1}{m-1} p_m^\frac{m}{m-1}$$
passing to the limit (using the a.e. convergence) yields $\rho_\infty p_\infty = p_\infty$ and thus
\begin{equation}
(\rho_\infty -1) p_\infty =0\quad\mbox{ a.e. } \RR_+\times \Omega,
\end{equation}
which gives the Hele-Shaw condition $p_\infty(x,t) \in P_\infty(\rho_\infty(x,t))$ a.e.  in $\RR_+\times \Omega$.

\medskip

Similarly, we have 
$$
\rho_m^m = \left(\frac{m-1}{m} p_m \right)^{\frac{m}{m-1}} \to p_\infty \quad\mbox{ a.e. } \RR_+\times \Omega.
$$
Since $\rho_m^m$ is bounded in $\mathrm{BV}(\RR_+\times \Omega_k)$, the convergence holds in $L^1_{loc}$ as well.  Rewriting \eqref{eq:1} as
$$
\pa_t \rho_m = \Delta p_m + \lambda \rho_m
$$
and passing to the limit, we  deduce
$$
\pa_t \rho_\infty = \Delta p_\infty + \lambda \rho_\infty \qquad \mbox{ in } \mathcal D'(\RR_+\times \Omega).
$$

\subsection{Bounds on the gradient of $p_m$ and convergence of $\rho_m$}
\begin{lemma}\label{lem:2.4}
There exists a constant $C$ independent  of $m$ such that
\begin{equation}\label{eq:gradp}
\int \int_{Q_T} |\na p_m|^2\, dx \, dt\leq CT.
\end{equation}
Furthermore, $\{\rho_m\}_{m\in \NN}$ is relatively compact in $C^s(0,T;H^{-1}(\Omega))$ for all $s\in(0,1/2)$.
\end{lemma}

\begin{proof}
Integrating the equation for the pressure
\eqref{eq:pressure} yields
$$
\frac{d}{dt} \int_{\Omega} p_m\, dx = -(m-2) \int_{\Omega} |\na p_m|^2\, dx + (m-1)\int_{\pa K} p_m \na p_m\cdot \nu dS + (m-1) \int_\Omega \lambda p_m\ dx.
$$
Using \eqref{eq:flux} and the fact that $p_m=\frac{m}{m-1}f$ on $\pa K$ we deduce
$$
 \int_{\Omega} |\na p_m|^2\, dx \leq -\frac{1}{m-2} \frac{d}{dt} \int_{\Omega} p_m\, dx + \frac{m-1}{m-2}C.
$$
Integrating in time and using \eqref{eq:L1p} we deduce \eqref{eq:gradp}. Using \eqref{eq:1}, we deduce that
$$ \pa_t \rho_m \mbox{ is bounded in } L^2(0,T;H^{-1}(\Omega)) .$$
Since $\rho_m$ is bounded in $L^\infty(0,T;L^1(\Omega)) $ and in $L^\infty(0,T;L^\infty(\Omega))$, we also have 
$$ \rho_m \mbox{ is bounded in } L^\infty(0,T;L^2(\Omega)). $$
Since $H^{-1}(\Omega)$ is compactly embedded in $L^2(\Omega)$, Lions-Aubin Lemma (see for example \cite{Lions,Amann})
implies 
that $\{ \rho_m\} $ is relatively compact in
$$ C^s(0,T; H^{-1}(\Omega)) \mbox{ for all $s\in(0,1/2)$}.$$
 \end{proof}
 
 Estimate \eqref{eq:gradp} implies in particular that $\na p_\infty(\cdot,t)\in L^2(\Omega)$ for a.e. $t>0$: it
will be useful in the proof of Theorem \ref{prop:obs}.

The compactness of $\{\rho_m\}_{m\in \NN}$ in $C^s(0,T;H^{-1}(\Omega))$ implies  $\rho_\infty \in C^s([0,\infty);H^{-1}(\Omega))$.
Furthermore, 
\eqref{eq:infty} and  \eqref{eq:L1rho} implies that $\rho_m(t)$ is bounded in $L^1(\Omega)\cap L^\infty(\Omega)$ and thus converges, up to a subsequence, weakly in $L^\infty(\Omega)$ to $\rho_\infty \in [0,1]$. We will see later in Section \eqref{sec:unique} that the limit density is unique (Proposition~\ref{prop:order}), which  implies that the whole original subsequence converge to $\rho_\infty$ weakly in $L^\infty(\Omega)$.

\section{Proof of Theorem \ref{prop:obs} and Proposition \ref{prop:comp}, \ref{prop:obstacle}}

In  this section, we use the notation $\rho_m$ and $p_m$ even though we are only considering convergent subsequences. Let us first introduce a lemma to be used in the proof of Theorem~\ref{prop:obs}.

\begin{lemma}\label{lem:pressureH1}
For all $t_0\geq 0$, $p^+(\cdot,t_0)\in H^1(\Omega)$ and 
$$
\int_{\Omega}  |\na p^+(x,t_0) |^2 \, dx  \leq\liminf_{\delta\to0} \frac 1 \delta  \int_{t_0}^{t_0+\delta} \int_{\Omega}   |\na p_\infty|^2  \, dx\, dt.
$$
\end{lemma}

\begin{proof}
First, using \eqref{eq:ineq1} with a nonnegative test function $v\in H^1(\Omega)$ supported in $U$ (see \eqref{eq:rho1}) and satisfying the boundary condition on $\pa K$, we get:
$$ \frac 1 \delta  \int_{t_0}^{t_0+\delta} \int_{\Omega} |\na p_\infty|^2\, dx\, dt \leq C,$$
for some constant $C$ independent of $\delta$ and so the $\liminf$ exists and is finite.
Given $T(x) \in \mathcal (D(\Omega))^n$, we can write
$$ 
\frac 1 \delta  \int_{t_0}^{t_0+\delta} \int_{\Omega}  p_\infty \div T\,dx\, dt  
= - \frac 1 \delta  \int_{t_0}^{t_0+\delta} \int_{\Omega} \na p_\infty \cdot T\, dx\, dt \leq \left(\frac 1 \delta  \int_{t_0}^{t_0+\delta} \int_{\Omega}  |\na p_\infty |^2\,dx\, dt \right)^{1/2}  \| T\|_{L^2(\Omega)},$$
and we can pass to the limit $\delta\to 0$, using \eqref{eq:tracep}, to get
$$ 
 \int_{\Omega}  p^+(x,t_0) \div T(x)\,dx 
\leq
 \left(\liminf_{\delta\to0}\frac 1 \delta  \int_{t_0}^{t_0+\delta} \int_{\Omega}  |\na p_\infty |^2\,dx\, dt \right)^{1/2}  \| T\|_{L^2(\Omega)},$$
 and the result follows.
\end{proof}

\begin{proof}[Proof of Theorem \ref{prop:obs}]
Given $t_0\geq0$ and  a function $v(x)$ in $E_{t_0}$, 
%we set $v_m := \frac{m}{m-1} v$ so that $p_m(x,t_0)=v_m(x)=\frac{m}{m-1} f(x,t_0)$ on $\pa K$.
and using the equation for the pressure \eqref{eq:pressure} and density \eqref{eq:1} we can write, in $\mathcal D'(\RR_+)$,
\begin{align}
& \int_{\Omega} \na p_m \cdot \na p_m - \rho_m \na p_m\cdot\na v-  \lambda(p_m-v)\, dx \nonumber \\
&\qquad \qquad = -\frac{1}{m-1}\left[ \frac{d}{dt}\int_\Omega p_m\, dx - \int_\Omega|\na p_m|^2\, dx\right]+ \frac{d}{dt} \int_\Omega v \rho_m\, dx -\int_\Omega \lambda(\cdot,t) v (\rho_m-1)\, dx\nonumber \\
& \qquad\qquad\quad + \int_{\pa \Omega} [p_m-\rho_m v ] \na p_m\cdot \nu\, dS.\label{eq:mainineq}
\end{align}
Formally at least, it is not difficult to see that the variational formulation of the obstacle problem \eqref{eq:var} follows by passing to the limit $m\to \infty$ and taking $t=t_0$ in \eqref{eq:mainineq}.
The rest of the proof is devoted to making this limit rigorous to derive \eqref{eq:var} (for all $t_0>0$).
\medskip

First, using the boundary condition, we note that the last term is equal to
$$
 \int_{\pa \Omega} \left[\frac{m}{m-1} f(x,t)-f(x,t)^{\frac{1}{m-1}} f(x,t_0)\right] \na p_m\cdot \nu\, dS
$$ 
and thus satisfies (using \eqref{eq:condbc}):
\begin{align}
\limsup_{m\to \infty} \left|  \int_{\pa \Omega} [p_m-\rho_m v] \na p_m\cdot \nu\, dS
\right| 
& \leq C \int_{\pa \Omega} |f(x,t)- f(x,t_0)|  \, dS\nonumber \\
& \leq C |t-t_0| . \label{eq:bdterm}
\end{align}

Next, it is clear that  our a priori estimates do not allow us to pass to the limit in \eqref{eq:mainineq} pointwise in time. So, given $\delta>0$, we integrate \eqref{eq:mainineq} with respect to $t\in (t_0,t_0+\delta)$ and pass to the limit $m\to \infty$. The left hand side of  \eqref{eq:mainineq}  satisfies
\begin{align*}
& \liminf_{m\to\infty}  \int_{t_0}^{t_0+\delta} \int_{\Omega} \na p_m \cdot \na p_m - \rho_m \na p_m\cdot \na v - \lambda(p_m-v)\, dx\, dt \\
& \qquad \geq \int_{t_0}^{t_0+\delta} \int_{\Omega} |\na p_\infty|^2 -  \na p_\infty \cdot \na v -  \lambda(p_\infty-v)\, dx\, dt,
\end{align*}
where we used in particular the fact that $\rho_m \na p_m$ is bounded in $L^2(Q_T)$ by \eqref{eq:gradp} and thus converges weakly to $\na p_\infty$ (since $\rho_m \na p_m = \na \rho_m^m$ converges to $\na p_\infty$ in $\mathcal D'(Q_T)$).
Using \eqref{eq:gradp} and \eqref{eq:L1p} to control the first term in the right hand side of \eqref{eq:mainineq}  and \eqref{eq:bdterm} for the last term, we deduce
\begin{align}
& \int_{t_0}^{t_0+\delta} \int_{\Omega} |\na p_\infty|^2 -  \na p_\infty \cdot\na v - \lambda(p_\infty-v)\, dx\, dt\nonumber \\
& \qquad\leq  \liminf_{m\to\infty}
 \int_\Omega v(x) [\rho_m(x,t_0+\delta) - \rho_m(x,t_0)] \, dx +  \int_{t_0}^{t_0+\delta}  \int_\Omega \lambda(\cdot,t) v(x) (1-\rho_\infty(x,t))\, dx dt + \mathcal O (\delta^2).\label{eq:ge}
%& \qquad\leq  \liminf_{m\to\infty}  \int_\Omega v(x) [\rho_m(x,t_0+\delta) - \rho_m(x,t_0)] \, dx+ \mathcal O (\delta^2)
\end{align}
%(since $v(x)\geq 0$ and $\rho_\infty(x,t) -1\leq 0$).

Formally, the first term in the right hand side is non-positive because $v \rho_\infty(\cdot,t_0+\delta) \leq v$ while $v\rho_\infty(\cdot,t_0)  = v$ (this is where we use the fact that $v\in E_{t_0}$). In order to make this rigorous, we first note that
$$ \frac{d}{dt}  \int_\Omega v(x) \rho_m(x,t)\, dx = - \int_\Omega \rho_m \na p_m \cdot \na v \, dx + \int_{\pa K} \rho_m v \na p_m\cdot \nu \, dS +\int_\Omega \lambda(\cdot,t) \rho_m v\, dx
$$
and so (using the fact that $\rho^{m-1}_m|_{\pa K} =v|_{\pa K} =f$):
$$\left|  \frac{d}{dt}  \int_\Omega v(x) \rho_m(x,t)\, dx\right| \leq \int_\Omega  |\na p_m(x,t)||\na v(x)| \, dx +  
\int_{\pa K}  f^{\frac{1}{m-1}}(x,t) f(x,t_0)|\na p_m\cdot \nu| dS +\Lambda  \int_\Omega \rho_m v(x)\, dx.
$$
The first term in the right hand side is bounded in $L^2(0,T)$ (using \eqref{eq:gradp}), and the second term is bounded in $L^\infty(0,T)$ (using \eqref{eq:flux}).
We deduce that the function
$t\mapsto \int_\Omega v_m(x) \rho_m(x,t)\, dx$ is bounded in $H^1(0,T)\subset C^{1/2}[0,T]$ and thus converges (up to a subsequence) uniformly in $[0,T]$. Since  $\int_\Omega v (x) \rho_m(x,t)\, dx$ converges to $ \int_\Omega v(x) \rho_\infty(x,t)\, dx$ in $\mathcal D'(\RR_+)$, we have
$$ \int_\Omega v(\cdot) \rho_m(\cdot,t)\, dx \to  \int_\Omega v(\cdot) \rho_\infty(\cdot,t)\, dx  \mbox{ locally uniformly in $\RR_+$}.$$
Consequently
$$
  \liminf_{m\to\infty}
 \int_\Omega v(x) [\rho_m(x,t_0+\delta) - \rho_m(x,t_0)] \, dx = \int_\Omega v(x) [\rho_\infty(x,t_0+\delta) - \rho_\infty(x,t_0)] \, dx \leq 0,
 $$
where we used the fact that $v(x) \rho_\infty(x,t_0) = v(x)$ (since $v\in E_{t_0}$) and $\rho_\infty\leq1$.

Going back to \eqref{eq:ge}, we deduce (using the fact that $v(x) (1-\rho_\infty(x,t))\geq 0$) : 
\begin{align}
& \int_{t_0}^{t_0+\delta} \int_{\Omega} |\na p_\infty|^2 -  \na p_\infty \cdot\na v -\lambda(p_\infty-v)\, dx\, dt
%& \qquad  \leq \int_\Omega v(x) [\rho_\infty(x,t_0+\delta) - \rho_\infty(x,t_0)] \, dx
%+ \int_{t_0}^{t_0+\delta}  \int_\Omega \lambda(\cdot,t) v(x) (\rho_\infty(x,t) -1)\, dx dt + \mathcal O (\delta^2)\nonumber \\
%& \qquad  \leq \int_\Omega v(x) (\rho_\infty(x,t_0+\delta)-1) \, dx+ \int_{t_0}^{t_0+\delta}  \int_\Omega \lambda(\cdot,t) v(x) (\rho_\infty(x,t) -1)\, dx dt+ \mathcal O (\delta^2)\nonumber \\
 \leq \Lambda  \int_{t_0}^{t_0+\delta}  \int_\Omega v(x) (1-\rho_\infty(x,t) )\, dx dt+ \mathcal O (\delta^2),\label{eq:ef}
\end{align}
To prove the result, it remains to divide by $\delta$ and pass to the limit $\delta\to 0$.
% using Lebesgue's differentiation theorem. To avoid having a set of valid $t$ that depends on $v$, 
We first use Young's inequality to rewrite \eqref{eq:ef} as:
\begin{align}
& \frac 1 \delta  \int_{t_0}^{t_0+\delta} \int_{\Omega}\frac 1 2  |\na p_\infty|^2 - \lambda p_\infty \, dx\, dt\nonumber\\
& \qquad \leq 
   \frac 1 \delta \int_{t_0}^{t_0+\delta} \int_{\Omega}  \frac 1 2  |\na v|^2 -  \lambda(\cdot,t) v\, dx\, dt +
    \frac \Lambda \delta \int_{t_0}^{t_0+\delta}  \int_\Omega   v(x) (1-\rho_\infty(x,t) )\, dx dt+ \mathcal O (\delta) \nonumber \\
 &\qquad  \leq  \int_{\Omega}  \frac 1 2  |\na v|^2 -  \int_{\Omega} \left(\frac 1 \delta \int_{t_0}^{t_0+\delta}\lambda(\cdot,t) dt \right) v\, dx + 
  \frac \Lambda \delta \int_{t_0}^{t_0+\delta}   \langle1- \rho_\infty(x,t), 
  v(x) \rangle_{H^{-1},H^1} dt +\mathcal O (\delta).\label{eq:ineq1}
\end{align}
Since  $\rho_\infty \in C([0,\infty);H^{-1}(\Omega))$ and $v\in E_{t_0}$, we have
$$ \lim_{\delta\to 0}   \frac 1 \delta \int_{t_0}^{t_0+\delta}   \langle1- \rho_\infty(x,t) , 
   v(x) \rangle_{H^{-1},H^1} dt  =   \langle 1- \rho_\infty(x,t_0) , 
   v(x) \rangle_{H^{-1},H^1} = 0 \quad \forall t_0\geq 0$$
and we can use \eqref{eq:tracep} (and a similar inequality for $\lambda$) to pass to the limit in the terms involving $\lambda$.

\medskip

From Lemma~\ref {lem:pressureH1}, we deduce that 
$$
\int_{\Omega}\frac 1 2  |\na p^+(x,t_0) |^2 - \lambda^+(x,t_0)  p^+(x,t_0) \, dx 
\leq \int_{\Omega}  \frac 1 2  |\na v(x)|^2 - \lambda^+(x,t_0) v(x)\, dx \qquad  \forall t_0>0,
$$
which implies that $p^+(\cdot,t)$ is indeed a solution of \eqref{eq:obstacle} for every $t\geq0$.

%Finally, Lebesgue differentiation theorem implies
%$$ \frac 1 \delta  \int_{t_0}^{t_0+\delta} \int_{\Omega}\frac 1 2  |\na p_\infty|^2 - \lambda p_\infty \, dx \longrightarrow \int_{\Omega}\frac 1 2  |\na p_\infty(x,t_0) |^2 - \lambda p_\infty(x,t_0) \, dx \qquad \hbox{ for a.e. } t_0>0.
%$$
%We deduce
%$$
%\int_{\Omega}\frac 1 2  |\na p_\infty(x,t_0) |^2 - \lambda p_\infty(x,t_0) \, dx 
%\leq \int_{\Omega}  \frac 1 2  |\na v|^2 - \int_{\Omega}\lambda(\cdot,t_0) v\, dx \qquad  \hbox{ for a.e. } t_0>0
%$$
%which implies that $p_\infty(\cdot,t)$ is indeed a solution of \eqref{eq:obstacle} for a.e. $t>0$.

\medskip

The derivation of \eqref{eq:var} is classical (given $u\in E_t$ and $\eps>0$, take $v = p+\eps(u-p)$ in  \eqref{eq:obstacle} and pass to the limit $\eps\to 0$). The uniqueness of $p^*$ follows from \eqref{eq:var}: if $p_1$ and $p_2$ are two solutions, then by plugging in each other as test functions we obtain
$$  \int_{\Omega}| \na (p_1-p_2)|^2 dx =0 ,$$
and thus $p_1=p_2$.
\end{proof}

\bigskip

\begin{proof}[Proof of Proposition \ref{prop:obstacle}]
For any ball $B_r(x_0)\in \mathcal O(t)$, we have  $\langle \vphi ,(1-\rho_\infty(t)\rangle =0$ for any $\vphi \in \mathcal D(B_r(x_0))$, and so $p$ solves the classical obstacle problem in $B_r(x_0)$. The usual theory (see \cite{caffarelli}) implies that $p\in C^{1,1}(B_{r/2}(x_0))$ and satisfies $\Delta p = \lambda \chi_{\{p>0\}}$ in $B_{r/2}(x_0)$. The proposition follows.

\end{proof}

\section{Uniqueness of the limit solution and comparison principle}\label{sec:unique}

In this section we establish the uniqueness for the limit problem in  a general bounded domain $D$ of $\RR^n$ (with smooth boundary):
Given  a  continuous  function $g\geq 0$  defined on $\partial D\times [0,T]$ and $\brho\geq 0$ a nonnegative function in $D$ satisfying  $0\leq \brho \leq 1$,  we consider the problem
\begin{equation}\label{eq:P}
\begin{cases}
\pa_t \rho  = \Delta p + \lambda \rho , & \mbox{ in } D\times (0,T], \qquad  p\in P_\infty(\rho) \mbox{ a.e.} \mbox{ in } D\times (0,T] ;\\
p =g &\hbox{ on }  \pa D\times [0,T];\\
\rho (t=0) = \brho & \mbox{ in } D.
\end{cases}
\end{equation}

A weak solution of \eqref{eq:P} is a set of functions  $(\rho,p)\in L^\infty(D\times (0,T]) \times L^2(0,T;H^1(D))$ satisfying \eqref{eq:P} in the sense of distribution.

In particular, the condition $ p\in P_\infty(\rho)$ implies that $0\leq\rho\leq 1$ and $p(1-\rho)=0$ a.e. in $D\times (0,T]$ and for any smooth, compactly supported test function $\psi: \Omega \to \RR$ with $\psi(\cdot,T)=0$ and $ \psi=0$ on $\partial D \times [0,T]$ we have
\begin{equation}\label{weak_limit}
\int_{D \times [0,T]} (\rho \psi_t + p \Delta\psi + \lambda \rho \psi) dxdt = -\int_D \rho(\cdot,0)\psi(\cdot,0) dx + \int_0^T\int_{\partial D} g\partial_{\nu}\psi dS dt.
\end{equation}

%\begin{definition}
% The {\it weak solution }$(\rho, p)$ of \eqref{eq:P}  is a pair of nonnegative, bounded functions in $D\times [0,T]$ that satisfies the following: 
%\begin{itemize}
%\item[(a)] The trace of $p(\cdot,t)$ equals  $f$ on $\partial D$ a.e. $t\in [0,T]$ \textcolor{black}{do we need to require $p\in H^1(D)$ to define the trace?};
%\item[(b)]  $0\leq \rho \leq 1$ and  $p(1-\rho)=0$ a.e. in $D\times (0,T]$;
%\item[(c)] For any smooth, compactly supported test function $\psi: \Omega \to \RR$ with $\psi(\cdot,T)=0$ and $ \psi(\cdot,0)=0$ on $\partial D \times [0,T]$ we have
%\begin{equation}\label{weak_limit}
%\int_{D \times [0,T]} (\rho \psi_t + p \Delta\psi - \lambda \rho \psi) dxdt = -\int_D \rho(\cdot,0)\psi(\cdot,0) dx +\int_0^T\int_{\partial D} f\partial_{\nu}\psi dS dt.
%\end{equation}
%\end{itemize}
%\end{definition}

We then have the following result, which implies in particular Proposition \ref{prop:uniqueness}:

\begin{proposition}\label{prop:order}
Suppose $\lambda \in L^2([0,T]; H^1(D))$, then there is at most one weak solution  $(\rho,p)$ of \eqref{eq:P}.

Furthermore, if $(\rho_i,p_i)$ for $i=1,2$ are two pairs of weak solutions of \eqref{eq:P} with boundary data $g_i$ and initial data $\brho_i$ and if $\brho_1 \leq \brho_2$ in $D$ and $g_1 \leq g_2$ on $\partial D \times [0,T]$, then $\rho_1 \leq \rho_2$ a.e. in $D\times [0,T]$.  
\end{proposition}

\begin{proof}
To show the uniqueness we follow the Hilbert dual argument developed in \cite{PQV}. Since the proof is largely parallel, we will only remark on necessary modifications due to the presence of the fixed boundary $\partial  K$.

\medskip

Suppose $\psi$ is a nonnegative test function. Let us denote $D_T:= D\times (0,T)$. Taking the differences of the weak formulation \eqref{weak_limit} for $(\rho_i,p_i)$ for $i=1,2$, we have 
\begin{align*}
&\int\int_{D_T} [(\rho_1-\rho_2) \pa_t \psi + (p_1-p_2)\Delta \psi +\lambda(\rho_1-\rho_2)] \\  
&\quad = -\int_D (\brho_1-\brho_2)(x)\psi(x,0) dx + \int_0^T\int_{\partial D} (g_1-g_2)\partial_{\nu}\psi dS dt \\
&\quad\geq  \int_0^T\int_{\partial D} (g_1-g_2)\partial_{\nu}\psi dS dt.
\end{align*}
Thus
\begin{equation}\label{eq:uni}
\int\int_{D_T} (\rho_1-\rho_2+p_1-p_2)[A\pa_t \psi + B\Delta\psi + \lambda A\psi] dxdt \geq  \int_0^T\int_{\partial D} (g_1-g_2)\partial_{\nu}\psi dS dt.\end{equation}
where $\nu$ denotes the outward normal at $\partial D$ and 
$$ 
 A = \frac{\rho_1-\rho_2}{\rho_1-\rho_2+p_1-p_2} \;, \qquad B = \frac{p_1-p_2}{\rho_1-\rho_2+p_1-p_2} .
$$
 As in \cite{PQV} we define $A=0$ whenever $\rho_1=\rho_2$ (even when $p_1=p_2$) and $B=0$ when $p_1=p_2$ (even when $\rho_1=\rho_2$). Note that $A, B\in [0,1]$ due to the fact that $\rho(1-p)=0$.
 
 \medskip
 
 Let  now $G$ be a compactly supported and nonnegative smooth function in $D \times [0,T]$. As in \cite{PQV} the idea is to solve  the dual problem
\begin{equation}\label{eq:dual}
\left\{\begin{array}{lll}
A\pa_t \psi + B\Delta \psi + \lambda \psi = -AG &\hbox{ in }& D\times [0,T);\\
\psi = 0  &\hbox{ on } & \pa D\times [0,T]; \\
\psi(\cdot,T) = 0. &\hbox{ in }& D.
\end{array}\right.
\end{equation}
If $A$ and $B$ were strictly positive, by backward-in-time maximum principle, one can verify that $\psi$ is nonnegative. Thus it follows that $\partial_\nu \psi \leq 0$ on $\partial D \times [0,T]$. Thus going back to \eqref{eq:uni} and using the fact that  $g_1 \leq g_2$, it follows that 
\begin{equation}\label{eq:G}
\int\int_{D_T} (\rho_1-\rho_2) (-AG) \geq  0.
\end{equation}
Since $G$ is arbitrary nonnegative smooth function, we conclude that $\rho_1\leq \rho_2$ a.e. in $D\times [0,T]$. 

\medskip

However, $A$ and $B$ can be degenerate, so the argument requires the approximation of the dual problem \eqref{eq:dual}, by a regularized uniformly parabolic, Dirichlet boundary value problem (see \cite{PQV} for detailed description of this approximation). 
As in \cite{PQV},  we then pass to the limit in the regularization to deduce \eqref{eq:G}.
The assumption $\lambda \in L^2([0,T]; H^1(D))$  is necessary to ensure that the regularized problem produces small errors.

\medskip

To show uniqueness, suppose that $(\rho_i,p_i)$ are two solutions of \eqref{eq:P} with the same boundary condition $f$ and initial condition $g$. Then $\rho_1=\rho_2$ follows from the density ordering property obtained above. Once we have this, the difference of the weak equations yield
$$
\int\int_{D_t} (p_1-p_2)\Delta \psi dx dt = 0.
$$
Now as in \cite{PQV} we can choose $\psi$ to approximate $p_1-p_2$ to conclude that $p_1=p_2$ a.e. in $D_T$.

\end{proof}

\begin{remark}
It is not immediately clear that the pressure satisfy the ordering property (i.e.  $p_1 \leq p_2$  in Proposition \ref{prop:order}). However, the characterization of the pressure given in Proposition \ref{prop:obstacle} implies that the  pressure ordering follows from the density ordering.
\end{remark}

Now let us state two consequences of this proposition, based on the comparison principle for \eqref{eq:1}. First let us discuss our original problem with $\Omega := \RR^n\setminus K$. Recall that from Lemma~\ref{lem:basic} that the support of $\rho_m$ lies in $B_{\overline R+ C(T)}$  for given time range $0\leq t\leq T$. Therefore, setting $R(T):= \overline R+ C(T)$, their limit solution $(\rho_{\infty}, p_{\infty})$ is a weak solution of \eqref{eq:P} with $D:= B_{R(T)}\setminus K$, $g=f$ on $\pa K$ and $g=0$ on $\pa B_{R(T)}$.  Therefore we have the following corollary:

\begin{corollary}\label{cor:first}
Given $T>0$,  any weak solution of \eqref{eq:P} with $D:= B_{R(T)}\setminus K$, $g=f$ on $\pa K$, $g=0$ on $\pa B_{R(T)}$ and initial data $\brho=\rho^0$  is the $L^1(Q_T)$- limit of the functions $(\rho_m, p_m)$ solutions of \eqref{eq:1}. 
In particular, it follows  that  the pressure ordering property is true in this setting. 
\end{corollary}

The next observation will be useful, when we construct radial limit solutions with explicit free boundary motion laws. 
\begin{corollary}[Comparison Principle]\label{cor:cp}
Let $(\rho_\infty,p_\infty)$ be the limit solution of \eqref{eq:1} in $\Omega \times [0,T]$.
 If $D$ is a domain with smooth boundary that does not  intersect $K$ and if $(\rho_1, p_1)$  is a weak solution of \eqref{eq:P} in $D\times [t_1,t_2]$., then the following holds: If $p_\infty \leq p_1  $  on $\partial D \times [t_1,t_2]$ and $\rho_\infty \leq \rho_1$ on $t=t_1$, then $p_\infty \leq p_1$ and $\rho_\infty \leq \rho_1$ in $D \times [t_1, t_2]$.
\end{corollary}

\begin{proof}
Since $D$ does not intersect  $K$,  it is easy to check that $(\rho_\infty,p_\infty)$ is a weak solution of \eqref{eq:P} in $D\times [t_1,t_2]$ with initial data $\rho_\infty(\cdot,t_1)$ and fixed boundary data given as the trace of $p_\infty$ on $\partial D \times [t_1, t_2]$ (such trace exists a.e. in time since $p_\infty(\cdot ,t)$ is in $H^1(D)$ a.e. $t>0$). Now we can conclude from  Proposition~\ref{prop:order}.

\end{proof}

\section{Proof of Proposition \ref{prop:mu}}
{In the sequel, we write $p(t)$ instead of $p_\infty(t)$ for the unique solution of the obstacle problem \eqref{eq:obstacle}. We also recall that $\mathcal P(t)=\{p(t)>0\}$}.
We first show that $\supp \mu_t \subset \pa \mathcal{P}(t)\setminus \mathcal{O}(t)$:
For all smooth test functions $\vphi\in \mathcal D(\Omega)$, by definition of $\mu_t$ we have
$$ \mu_t(\vphi) = \int_\Omega (-\na p \cdot \na \vphi +\lambda(\cdot,t) \chi_{\mathcal{P}(t)} \vphi)\, dx.$$
Clearly, if $\vphi$ is supported in $\{ p(\cdot,t)=0\}$, the fact that $p\in H^1(\Omega)$ implies that $\na p=0$ a.e. in $\{p=0\}$ and thus $ \mu_t(\vphi)=0$.
And if $\vphi$ is supported in $\mathcal O(t)$, \eqref{eq:pob} implies
$$ \mu_t(\vphi) =0.$$
Since $\mathcal O(t)$ is an open set, we deduce that
$$ \supp(\mu_t) \cap \mathrm{Int} (\{ p(t)=0\}) =\emptyset, \qquad \supp(\mu_t) \cap \mathcal O(t)=\emptyset.
$$

On the  other hand note that $\mathrm{Int}(\mathcal{P}(t))\subset \mathcal O(t)$. Indeed if $p(t)>0$ in $B_\delta(x_0)$, then $1-\rho_\infty(t)=0$ a.e. in $B_\delta(x_0)$ and so $\int_{B_\delta(x_0)} (1-\rho_\infty(t))\, dx=0$. It follows that $x_0\in \mathcal O(t)$. Thus we can conclude that 
 $\mu_t$ is supported in  $\pa\mathcal{P}(t)\setminus \mathcal O(t)$ as claimed in Proposition \ref{prop:mu}.

\medskip

Next we show that $\mu_t$ is nonnegative. 
\begin{comment}
In dimension $2$ and $3$, the continuity of $p(\cdot,t)$ for a.e. $t>0$ and Proposition \ref{prop:obstacle} imply that $p$ solves $\Delta p(\cdot,t)  =\lambda(\cdot,t)$ in the open set $\mathcal{P}(t)$ for a.e. $t>0$.
Together with \eqref{eq:Dp}, this implies that $\Delta p(\cdot,t)  \geq \lambda\chi_{\mathcal{P}(t)}$ and therefore $\mu_t\geq 0$ (a.e. $t>0$).
In higher dimensions we can prove that $\mu_t\geq 0$ as follows:
\end{comment}
Define the function
$$Q_\delta(s):=\begin{cases}
\frac{s}{\delta} & \mbox{ if } s\in[0,\delta];\\
1 & \mbox{ if } s\geq \delta.
\end{cases}
$$
For any test function $\vphi \in \mathcal D(\Omega)$ satisfying $0\leq \vphi(x)\leq 1$, we write
\begin{align*}
\mu_t(\vphi) 
& = \int_\Omega - \na p \cdot \na \vphi + \lambda \chi_{\mathcal{P}(t)} \vphi\, dx\\
& = \int_\Omega -\na p \cdot \na (\vphi Q_\delta(p)) + \lambda \vphi Q_\delta(p) \, dx
 + \langle \Delta p, \vphi (1-Q_\delta(p)) \rangle  + \int_\Omega \lambda \chi_{\mathcal{P}(t)}\vphi(1-Q_\delta(p)) \, dx.
\end{align*}
Using \eqref{eq:var} with $u= p-\delta \vphi Q_\delta(p)$ (which satisfies $p\geq u \geq p(1-\vphi)\geq 0$ and is thus admissible) the first integral is non-negative. Next note that
 $$
 \langle \Delta p, \vphi (1-Q_\delta(p)) \rangle  =  \int (\nabla p \cdot \nabla  \varphi(Q_\delta(p)-1) + \nabla p \cdot \varphi Q'_{\delta}(p) \nabla p) dx.
 $$
 The second term in above equality is nonnegative since $Q_{\delta}$ is increasing.
 For the first term, we note that $ \nabla  \varphi(Q_\delta(p)-1)$ converges a.e. to $\nabla  \varphi \chi_{\{p=0\}}$.  Lebesgue dominated convergence theorem implies that it converges in $L^2$ and thus the first term converges to
zero since $\na p =0 $ a.e. in $\{ p=0\}$. 
 
 \medskip

  Thus
$$
\mu_t(\vphi) \geq \int_\Omega \lambda \chi_{\mathcal{P}(t)}\vphi(1-Q_\delta(p)) \, dx.$$
Finally, we have $ \chi_{\mathcal{P}(t)} (1-Q_\delta(p)) \to 0$ a.e. in $\Omega$ when $\delta\to0$. Sending $\delta\to 0$ and using Lebesgue dominated convergence theorem, we can conclude
that $\mu_t(\vphi) \geq 0$ and the result follows.

%%%%%%%%%%%%%%%%%%%%%
%%%%%%%%%%%%%%%%%%%%%

\section{The velocity law} \label{sec:velocity}

\medskip

 In this section we determine the velocity law of the congested zone $\{\rho_\infty=1\}$ for the limit solution $(\rho_\infty,p_\infty)$ by using comparison principle and barriers, as in the usual viscosity solutions approach.
 First we will define the relevant notion of barriers   and  
 prove that the usual comparison with barriers holds for our limit solution $(\rho_\infty,p_\infty)$
(see Corollary \ref{cor:sub}-\ref{cor:sup}).
In Section \ref{sec:radial} we show that in the radial symmetric case, the barriers we construct are indeed classical solutions.

\subsection{Comparison with barriers} 

The difficulty in making \eqref{eq:twoeq} rigorous is the lack of regularity of the pressure or density interface ($\pa \mathcal{P}$ or $\pa \Sigma$) and the lack of monotonicity of its motion.
In this section, we construct sub- and super-solutions of the limiting problem \eqref{eq:1} to be used as barrier in a viscosity solution type approach.

\medskip

Let $ B_r$ be a ball  in $\Omega$, and let $D$ be either $\Omega\setminus \overline B_r$ or $ B_r$.  For  a given time interval $[t_1,t_2] \subset [0,\infty)$ we consider a function (the pressure) $\phi \in C_c(\overline {D}\times [t_1,t_2])$ such that $\{\phi(t)>0\}$ is monotone (increasing or decreasing) and an initial density $\rho_1(x)$ satisfying $\rho_1 = 1$ in $\{\phi(t_1)>0\}$. 
We assume that $\{\phi(t)>0\}$ and $\rho_1(x)$ are such that the external density $\rho^E_\phi$, defined below, satisfies  
\begin{equation} \label{eq:rhoEcond}\rho^E_\phi (x,t)<1\mbox{ in } \{\phi=0\}.
\end{equation}

\medskip

This external density $\rho^E_\phi(x,t)$ solves the equation $\pa_t \rho =  \lambda\rho$ in the (decongestion) set $\{\phi=0\}$ together with appropriate boundary conditions. This leads to the following definitions:

If $\{\phi(t)>0\}$ is {\bf increasing} (``expanding solution"), then for all $x\notin \{ \phi(t_1)>0\}$, we define $t(x) = $ the last time that $\phi(x,t)=0$ (with $t(x)=t_2$ is $\phi(x,t_2)=0$) and set
$$ \rho^E_\phi(x,t) = \rho_1(x) \exp\left( \int_{t_1}^t \lambda (x,s) \, ds\right) \qquad \mbox{ for all $t<t(x)$}$$
(condition \eqref{eq:rhoEcond} is satisfied if $\rho_1(x)$ is small enough in $\{\phi(t_1)=0\}$).

If $\{\phi(t)>0\}$ is {\bf decreasing} (``contracting solution"), then for all $x\notin \{ \phi(t_2)>0\}$, we define $t(x)=$ the first time that $\phi(x,t)=0$ (with $t(x)=t_1$ is $\phi(x,t_1)=0$) and set
$$ \rho^E_\phi(x,t) = \rho_1(x) \exp\left( \int_{t(x)}^t \lambda (x,s) \, ds\right) 
\qquad \mbox{ for all $t>t(x)$}.$$
(condition \eqref{eq:rhoEcond} requires $\rho_1(x)$ to be small enough in $\{\phi(t_1)=0\}$, but since 
 $\rho_1=1$ in $\{\phi(t_1)>0\}$, it also requires $\exp\left( \int_{t(x)}^t \lambda (x,s) \, ds\right)<1$ for $x\in \{\phi(t_1)>0\} $).

\medskip

In both cases, we define the density in $D\times (t_1,t_2)$ by
\begin{equation}\label{eq:densitybarrier}
\rho_\phi (x,t) := \chi_{\{\phi(t)>0\}}(x) + \rho^E_\phi (x,t) (1-\chi_{\{\phi(t)>0\}}(x) ) = 
\begin{cases}
1 & \mbox{ in } \{\phi>0\}\\
 \rho^E_\phi (x,t)  & \mbox{ in } \{\phi=0\}.
 \end{cases}
\end{equation}

We then have:  
\begin{proposition}\label{prop:sub}
With the notation above,
assume that $(\rho_\phi,\phi)$  are such that 
\item[(a)]   $\phi\in C^1(\overline{\{\phi>0\}})\cap C^2_{loc}(\{\phi>0\})$ and $\Gamma:= \partial\{\phi>0\}$ is $C^2$ in space and $C^1$ in time.
\item[(b)] $\phi$ satisfies
\begin{equation}\label{eq:cbarrier}
\begin{cases}
-\Delta \phi  \leq \lambda & \mbox{ in } \{\phi>0\}; \\
(1-\rho_\phi^E) V_\phi \leq  |\na \phi| & \mbox{ on } \pa \{\phi>0\},
\end{cases}
\end{equation}
where $V_\phi$ denotes the normal velocity of the interface $ \pa \{\phi>0\}$.
\item Then $(\rho_\phi,\phi)$ is a  weak {\em subsolution} of the limiting problem \eqref{eq:P} in  $D\times [t_1,t_2]$, namely
$$
\pa_t \rho_\phi  \leq \Delta \phi + \lambda \rho_\phi  \mbox{ in }  D\times (t_1,t_2), \qquad  \phi \in P_\infty(\rho_\phi) \mbox{ a.e.} \mbox{ in } D\times (t_1,t_2)
$$
where the first equation holds in the sense that for every smooth, compactly supported test function $\psi: D\times (t_1,t_2) \to \RR$ with $\psi(\cdot,t_2)=0$ and $ \psi(\cdot,t)=0$ on $\partial D \times [t_1,t_2]$ we have
\begin{equation}\label{weak}
\int_{D \times [t_1,t_2]} (\rho_\phi \psi_t + \phi \Delta\psi + \lambda \rho_\phi \psi )dxdt \geq  -\int_{D} \rho_1(x)\psi(\cdot,t_1) dx + \int_{t_1}^{t_2}\int_{\partial B} \phi \partial_{\nu}\psi dS dt.
\end{equation}
\end{proposition}
Similarly, we have 
\begin{proposition}\label{prop:sup}
With the notation above,
assume that $(\rho_\phi,\phi)$  are such that
\item[(a)]   $\{\phi(\cdot,t)>0\}\Subset \Omega$ for all $t$, $\phi\in C^1(\overline{\{\phi>0\}})\cap C^2_{loc}( \{\phi>0\})$ and the interface $\Gamma:= \partial\{\phi>0\}$ is $C^2$ in space and $C^1$ in time.
\item[(b)] $\phi$ satisfies
\begin{equation}\label{eq:sbarrier}
\begin{cases}
-\Delta \phi \geq \lambda & \mbox{ in } \{\phi>0\};\\
(1-\rho_\phi^E) V_\phi \geq  |\na \phi| & \mbox{ on } \pa \{\phi>0\}.
\end{cases}
\end{equation}

Then $(\rho_\phi,\phi)$ is a {\em supersolution} of the limiting problem \eqref{eq:P}  in $D\times [t_1,t_2]$, namely
$$
\pa_t \rho_\phi  \geq \Delta \phi + \lambda \rho_\phi  \mbox{ in } D\times (t_1,t_2), \qquad  \phi \in P_\infty(\rho_\phi) \mbox{ a.e.} \mbox{ in } D\times (t_1,t_2). 
$$
(with the corresponding weak formulation as in \eqref{weak})
\end{proposition}

Note that for the contracting barrier, we have $V_\phi\leq 0$ and $\rho^E_{\phi}=1$ on $ \pa \{\phi(t)>0\}$ and so the free boundary condition reduces to $|\na \phi| \geq 0$ for subsolution and $|\na \phi| =0$ for supersolution.

\begin{proof}[Proof of Proposition \ref{prop:sub}]
We denote $S(t):= \{\phi(\cdot,t)>0\}=\{\rho(\cdot,t)=1\}$ and $\Gamma(t) = \partial S(t)\cap D$.
We also denote $\nu$ as the outward normal of the boundary of either $\Gamma(t)$ or $\partial D$ with respect to the domain $S(t)$. With these notations, we have 
\begin{align*}
\int_{D} \phi \Delta \psi \, dx = \int_{S(t)} \phi \Delta \psi \, dx  & \geq -\int_{S(t)} \lambda \psi \, dx - \int_{\pa S(t)} \psi \na \phi\cdot \nu \,  dS + \int_{\partial B} \phi \partial_{\nu}\psi dS\\
&\geq  -\int_{S(t)} \lambda \psi \, dx + \int_{\Gamma(t)}\psi  |\na \phi|\, dS + \int_{\partial B}\phi  \partial_{\nu}\psi dS,
\end{align*}
where we used the fact that $\phi=0$ and $\na \phi = |\na \phi| \nu$ on $\Gamma(t)$.

Next 
\begin{align*}
\int_{D} \phi \psi_t \, dx &= \int_{S(t)} \psi_t \, dx + \int_{D\setminus S(t)} \rho^E_\phi \psi_t \, dx     \\ 
 &= \frac{d}{dt} \int_{S(t)}  \psi \, dx-\int_{\Gamma(t)} V_\phi \psi \, dS + \int_{D\setminus S(t)} (\rho^E_\phi \psi)_t \, dx   - \int_{D\setminus S(t)} (\rho^E_\phi)_t \psi\, dx\\
 &= \frac{d}{dt} \int_D \rho_\phi \psi \, dx-\int_{\Gamma(t)} V_\phi(1-\rho^E)\psi \, dS  - \int_{D\setminus S(t)} (\rho^E_\phi)_t \psi\, dx\\
 &\geq  \frac{d}{dt} \int_D \rho_\phi\psi \, dx-\int_{\Gamma(t)} |\na \phi| \psi \, dS  - \int_{D\setminus S(t)} (\rho^E)_t \psi\, dx
 \end{align*}
Using the fact that
\begin{equation}\label{ode1}
(\rho^E)_t = \lambda \rho^E \hbox{  in }\{\phi=0\},
\end{equation} 
and the definition of $\rho_\phi$, we deduce
$$
\int_D (\rho_\phi \psi_t + \phi\Delta\psi) dx \geq -\int_D \lambda \rho_\phi \psi dx +\frac{d}{dt} \int_D \rho\psi dx + \int_{\partial B} \partial_{\nu}\psi dS,
$$
and we conclude by integrating with respect to $t\in(t_1,t_2)$.
\end{proof}
 
 The proof of Proposition \ref{prop:sup} is parallel.
  Note that it is not necessary to work with barriers such that the set $\{\phi(\cdot,t)>0\}$ is monotone: we chose to do so because the definition of $\rho^E_\phi$
 is more manageable in that case.
\medskip

Combining Proposition \ref{prop:sub} with the comparison principle for weak solutions of the limiting problem (Corollary~\ref{cor:cp}) we get: 

\begin{corollary}\label{cor:sub}
Let $(\rho_\phi,\phi)$ be as in Proposition \ref{prop:sub} (sub-solution). If 
\item[(i)] $\rho_1\leq \rho(\cdot,t_1)$ in $D$,
(so in particular $\{\phi(\cdot,t_1)>0\}\subset \{\rho_\infty(\cdot,t_1) =1\}$);
\item[(ii)] $\phi\leq p_\infty$ on $\pa B_r\times[t_1,t_2]$,
\item then $\rho_\phi \leq \rho_\infty$ in $D  \times [t_1,t_2]$.
In particular 
$$\{\phi(\cdot,t)>0 \} \subset \{\rho_\infty (\cdot,t)=1\}\mbox{ for all } t\in [t_1,t_2].
$$
\end{corollary}
Formally, this corollary says that a classical subsolution of the viscosity law (satisfying \eqref{eq:cbarrier}) cannot touch $\rho_\infty$ from below. In other words, $\rho_\infty$  satisfies the motion law
$$  (1-\rho^E_\infty)V_\infty \geq |\na p_\infty| \mbox{ in a viscosity sense.}$$

\begin{comment}
\textcolor{black}{is the  paragraph below to explain what we mean "in a viscosity sense?"}
We note that if a barrier is such that $\pa \{\phi(\cdot,t)>0 \}$ touches the  interface $\pa \{\rho_\infty (\cdot,t)=1\}$ from inside the set $\{\rho_\infty (t)=1\}$, then at the contact point, we have (assuming everything is smooth)
$ V_\infty \leq V_\phi $, $|\na p_\infty|\geq |\na \phi|$ and $\rho^E_\infty\geq \rho^E_\phi$ and so
$$  |\na p| -(1-\rho^E_\infty )V_\infty \geq   |\na \phi|-(1-\rho^E_\phi) V_\phi.$$
\end{comment}

\medskip

Similarly, Proposition \ref{prop:sup} implies:
\begin{corollary}\label{cor:sup}
Let $(\rho_\phi,\phi)$ be as in Proposition \ref{prop:sup} (super-solution). If 
\item[(i)] $\rho_1\geq \rho(\cdot,t_1)$ in $D$ (so in particular $\{\rho_\infty(\cdot,t_1) =1\} \subset \{\phi(\cdot,t_1)>0\}$)
\item[(ii)] $\phi\geq p_\infty$ on $\pa B_r\times[t_1,t_2]$
\item Then $\rho_\phi \geq \rho_\infty$ in $D  \times [t_1,t_2]$.
In particular 
$$
\{\phi(\cdot,t)>0 \} \supset \{\rho_\infty (\cdot,t)=1\}\mbox{ for all } t\in [t_1,t_2].
$$
\end{corollary}
As above, this result should be interpreted as saying that $\rho_\infty$ satisfies
$$  (1-\rho^E_\infty)V_\infty \leq |\na p_\infty| \mbox{ in a viscosity sense.}$$

 %\begin{corollary}\label{cp:2}
%Let $(\rho,p)$ be a weak solution of \eqref{eq:P} in $\Omega \times [0,T]$, and let $\phi$ be an expanding or contracting barrier in $D_T:=D\times [0,T]$, where $D$ is a open subset of $\Omega$ with a smooth boundary. Then $\rho$ cannot cross $\rho_{\phi}$ given in \eqref{associated} from below or above in $D_T$, if they are ordered on the parabolic boundary of $D_T$.  
%\end{corollary}'

Typically, for free boundary problems such "barrier property" allows us to introduce a notion of viscosity solutions which describes the pointwise behavior of the interface via comparison with barriers (see e.g.\cite{KP}).  It is thus natural to ask whether our weak solutions coincide with viscosity solutions. While we suspect that viscosity solutions theory can be  established for our problem, answering this question would require a different set-up of function spaces, and we do not pursue this question here to keep our investigation focused.

\subsection{The radial symmetric case}\label{sec:radial}

In this section we show that the free boundary velocity law holds in the classical sense in the radial setting as long as $\pa_t \lambda$ does not change signs too often. To simplify our discussion we further assume that $\lambda$ is  non-positive, since construction of radial barriers for positive $\lambda$ has been carried out in \cite{KP}.

\medskip

We thus assume that $K=B_1$ and that the boundary data is constant (we can take $f=1$ without loss of generality) and for simplicity we take  $\lambda= \lambda(t)\leq 0$ independent of $x$ monotone $C^1$ function of $t$. 
The analysis could be extended to radial symmetric functions $\lambda(|x|,t) \leq 0$ such that $\pa_t \lambda$ changes sign a finite number of time in the interval $[0,T]$.
\medskip

In this setting, we  construct compactly supported, radial  sub and super solutions of \eqref{eq:2} in $Q_T:= \{ |x| \geq 1\}\times [0,T]$.

For a given $R>1$, let us define $\phi_R(\cdot,t)$ as a solution of the Dirichlet boundary problem in $1\leq |x|\leq R$:
\begin{equation}\label{eq:phiR}
-\Delta\phi = \lambda(t) \hbox{ in } |x|<R,  \quad \phi=0 \hbox{ on } |x|=R, \quad \hbox{ and }\phi=1\hbox{ on }|x|=1.
\end{equation}
Note that this function will take negative value if $R$ is large (depending on $\lambda$).

For a given $R_0>1$, we assume that the initial density $\rho_0$ equals $1$ on $1\leq |x|<R_0$ and is strictly less than $1$ and Lipschitz in $|x|\geq R_0$. We assume that $R_0$ is small enough so that 
the initial pressure $\phi_{R_0}(\cdot,0)$ is nonnegative.
We then define the external density in the region $|x|\geq R$ by
$$
\rho^E(|x|,t):= \rho_0(|x|) \exp\left({\int_0^t \lambda(s) ds}\right) <  1 \hbox{ in } |x| \geq R_0,
$$
then $\rho^E(\cdot,t)$ is Lipschitz continuous. It is also straightforward to check that the function $\partial_r \phi_R(R)$ is Lipschitz continuous for $ R_0 <R < \infty$. Thus we can solve the following ODE for $0\leq t\leq T$:
\begin{equation}\label{velocity_law}
R'(t) = F(R(t),t), \hbox{ where } F(R,t):= \dfrac{(\partial_r\phi_R)_-(R,t)}{1-\rho^E(R,t)}, \qquad R(0)=R_0.
\end{equation}

Note that $\partial_r\phi_R(\cdot,t) \geq 0$ if and only if the function $\phi_R(\cdot,t)$ has a negative minimum in $1\leq |x|\leq R$. Indeed if $\partial_r \phi_R(R,t)<0$ and $\phi_R(\cdot,t)$ takes negative minimum,  from the radial symmetry of $\phi_R(\cdot,t)$ it follows that the function has a local positive maximum  for some $x$ such that $1 <|x|< R$, which contradicts  the subharmonicity of $\phi_R(\cdot,t)$. 

So as long as $\phi_{R(t)} $ is a non-negative function, we have 
$\partial_r\phi_R(\cdot,t) < 0$ and $R'(t)>0$ (which provides is an expanding solution of the limiting problem. and we can show that it happened when the function $t\mapsto \lambda(t)$ is decreasing.

\medskip

\noindent{\it Case 1:} {\it  $t\to \lambda(t)$  is increasing}:
In this case, we can define
\begin{equation}\label{1}
\phi(\cdot,t):= \phi_{R(t)}(\cdot,t) \hbox{ for } 0\leq t\leq T
\end{equation} 
and  we claim that $\phi$ stays nonnegative for all times. 

To show this, suppose that $\phi(\cdot,t)$ has a negative minimum at some time $t=t_0$. Then by continuity of $R(t)$, the same is true for $\phi(\cdot,s)$ for $s$ sufficiently close to $t_0$.  Hence from above discussion we have $R'(t)=0$ in a small time interval $[t_0-\e,t_0]$. Suppose we choose $\e$ such that $\phi(\cdot,t_0-\e)$ no longer has negative minimum. This must be true at least with $\e=t_0$ due to our assumption. But since $\lambda(t_0-\e) < \lambda(t)$ and $R(t_0-\e) = R(t)$ for $s=t_0-\e$, we have $\phi(x,t_0-\e)  \leq \phi(x,t_0)$, which is a contradiction to our choice of $\e$. 

Hence we have shown our claim, and it follows from \eqref{velocity_law} and Propositions \ref{prop:sub}
-\ref{prop:sup} that $\phi$ is an expanding solution of \eqref{eq:2} for all $t\geq 0$.

\medskip

\noindent{\it Case 2:} {\it  $t\to \lambda(t)$  is non-increasing}: In this case, $\phi(\cdot,t)$ might take negative value for some positive time.
We thus define
$$
t^*:= \sup\{t\in [0,T]: \phi(\cdot,t) \geq 0 \hbox{ in } 1\leq |x| < R(t)\}.
$$
If $t^*=\infty$ then we can define $\phi$ by \eqref{1} as above. We thus assume that $t^*<\infty$.
The same arguments as above implies that $|D\phi|=0$ at $(R(t^*),t^*)$.
Since $\lambda$ is non-increasing, it follows that $\phi_{R(t^*)}(\cdot,t)$ turns negative for $t>t^*$. 
For $t \geq t^*$ we define $\tilde{R}(t)$ as the unique boundary point of $\{\psi(\cdot,t)>0\}$, where $\psi(\cdot,t)$ solves the obstacle problem 
$$
-\Delta \psi= \lambda(\cdot,t) \chi_{\{\psi>0\}} \hbox{ in } 1<|x|<R(t^*), \quad \hbox{ with }\psi =1\hbox{ on } |x|=1.
$$
  We then define
 
 \begin{equation}\label{2}
 \phi(\cdot,t) := \phi_{R(t)}(\cdot,t) \hbox{  for } 0\leq t\leq t^*, \quad  \phi(\cdot,t):= \phi_{\tilde{R}(t)} (\cdot,t) \hbox{ for } t^* \leq t \leq T.
 \end{equation}
 Since $\lambda$ is non-increasing,  so is $\tilde{R}$ and  $|D\phi|(\tilde{R}(t),t) = 0$. It follows that $\phi$ is a contracting solution for $t^* \leq t\leq T$. 
 
\medskip

 Below is the summary of our  conclusion:

\begin{lemma}
\item If $t\to \lambda(t)$ is increasing, then the function $\phi$  defined by \eqref{1} is an expanding solution for $0\leq t\leq T$. 
\item If $t\to\lambda(t)$ is non-increasing, then the function $\phi$ defined in \eqref{2}   is an expanding solution for $0\leq t\leq t^*$ and is a contracting solution for $t^*\leq t\leq T$. 
\end{lemma}

 Due to the uniqueness of the limit problem we can now completely characterize the limiting profile of radial solutions for $\lambda$ that are monotone $C^1$ function of time.
 
\begin{proposition}\label{existence}
Assume that $K=B_1$, $f=1$ and that $t\mapsto \lambda(t)$ is  a monotone $C^1$ function. 
Let $\rho_0^m$ be a radially symmetric function satisfying the conditions of Assumption \ref{ass:init}.
Then the limit $(\rho_\infty,p_\infty)$  given by Theorem \ref{thm:1} is radially symmetric and satisfies
$$
\begin{cases}
\Delta p_\infty + \lambda = 0 & \mbox{ in } \{p_\infty>0\}; \\
(1-\rho_\infty^E) V  \leq  |\na p_\infty| & \mbox{ on } \pa \{p_\infty>0\}.
\end{cases}
$$
Furthermore
\begin{itemize}
\item[(a)] If $t\mapsto \lambda(t)$ is increasing, then $\{p_\infty>0\}$ is always expanding ($\rho_\infty^E<1$,  $|\na p_\infty| >0$ and $V>0$ on $\pa \{p_\infty>0\}$)
\item[(b)] If $t\mapsto \lambda(t)$ is non-increasing, then there exists a time $t^*\in [0,T]$ such that $\{p_\infty>0\}$ is expanding for $0\leq t\leq t^*$ and contracting for $t^* \leq t\leq T$
($\rho_\infty^E=1$,  $|\na p_\infty| =0$ and $V\leq 0$ on $\pa \{p_\infty>0\}$)
\end{itemize}
\end{proposition}

\subsection{Continuous expansion of the congested zone}

As an application of comparison principle (Corollary~\ref{cor:cp}) with a radial barrier, we show that the congested zone does not expand discontinuously over time. Note that it may shrink discontinuously even if $\lambda$ is smooth, for instance due to topological changes. Note also that if $\lambda$ is nonnegative, the expansion may not be continuous due to the nucleation of congested zones created by the growth of external densities.

\begin{corollary}
If $(\rho,p)$ is a limit solution in $\Omega \times [0,T]$ and $\lambda\in C(Q_T)\cap L^2([0,T]; H^1(\Omega))$ is negative, then
\begin{equation}\label{no_jump}
\overline{\{p>0\}\cap Q_T} =\overline{\{p>0\}} \cap Q_T = \overline{\{\rho=1\}}\cap Q_T\hbox{ for any } T>0.
\end{equation}
\end{corollary} 

\begin{proof}
We denote
$$S_1:=\overline{\{p>0\}\cap Q_T}, \qquad S_2:=\overline{\{p>0\}}\cap Q_T.$$
Since $S_1\subset S_2$ by definition, we only need to show that $S_2 \subset S_1$ in order to prove the first equality.

Given $x_0\notin S_1$ there exists $r>0$ such that
$$
B_{2r}(x_0) \times [T-r, T) \subset \{p=0\}.
$$
We claim that $B_{r/2}(x_0)$ lies in $\{p(\cdot,T)=0\}$. 
This proves that $(x_0,T)\notin S_2$, hence $S_2 \subset S_1$.

\medskip

To show that $B_r(x_0) \subset \{p(\cdot,T)=0\}$, we use a barrier argument in $\Sigma:=B_r(x_0)\times [T-\e, T)$ for a sufficiently small $\e>0$ as follows.  Due to Proposition~\ref{prop:mu} we have $\rho_t = \lambda \rho$ in $B_{2r}(x_0) \times [T-r,T)$, and thus 
$$
\rho <a(\lambda, r)<1 \hbox{ in }B_{2r}(x_0) \times [T-r/2, T).
$$  Let us construct an expanding supersolution in $\Sigma$ as follows. Let  $\phi_0$ solve
$$
-\Delta\phi_0= \Lambda\hbox{ in } \{r<|x|<2r\},  \quad \phi_0=0 \hbox{ on } \{|x|=r\}, \quad \hbox{ and }\phi_0=M:= \|p\|_{L^{\infty}(\Omega \times [0,T])} \hbox{ on } \{|x|=2r\},
$$
and let $\phi(\cdot,t):=\phi_{R(t)}$ defined by \eqref{eq:phiR} with $\lambda =\Lambda$ where $R(t)$ solves
$$
R'(t)  = \frac{ |D\phi|(R(t),t)|}{1-a\exp^{(\Lambda t)}} \quad \hbox{ for } 0\leq t\leq \e,  \quad  \mbox{ with } R(0) = r.
$$
Then $\phi$ is an expanding supersolution in $\Sigma$ with fixed boundary data $M$  on $\partial B_{2r}(x_0)$ and initial data $\rho_0 = \chi_{r<|x|<2r} + a\chi_{|x| \leq r}$. Corollary~\ref{cor:cp} now applies to show that $p(\cdot,T)\leq \phi(\cdot,T)$. Choosing $\e =\e(M,a)$ sufficiently small so that $R(T) \leq R(0) + \frac{r}{2}$, it follows that $\phi(\cdot,T) = 0$ in $B_{r/2}(x_0)$ and we can conclude.

\medskip

It remains to show the second equality of the Corollary. Note that we have $\{p>0\} \subset \{\rho=1\}$, and thus their closures are also ordered. On the other hand we showed above that if $x_0$ lies outside of $\overline{\{p>0\}}$ then $\rho$ is strictly less than one in a small neighborhood of $x_0$, and thus it is outside of $\overline{\{\rho=1\}}$. The result follows.

\end{proof}

\section{Monotone increasing solutions}

In this section  we suppose that $\lambda \in L^2([0,T];H^1(\Omega))$ is non-decreasing in time.  We first show that in this setting, if the density starts  as a characteristic function, the pressure only increases over time.

\begin{lemma}\label{patch}
Let $\Sigma_0$ be a bounded subset of $\RR^n$ which contains $K$. Suppose that $\rho_0= \chi_{\Sigma_0\setminus K}$ and that $\Sigma_0\setminus K$ coincides with the initial pressure support  $\{p_0>0\}$, where $p_0$ solves  \eqref{eq:var} with $\rho_{\infty}(\cdot,t)$ replaced by $\rho_0$. If $(\rho,p)$ is the  limit solution given by Theorem~\ref{thm:1} with initial data $\rho_0$, then $\rho$ and $p$ are monotone  increasing with respect to $t$.
\end{lemma}

\begin{proof}
Let $B_R$ contain the support of $\Sigma_0$. We claim that $(\rho_0,p_0)$ is a stationary subsolution of \eqref{eq:P} with $D = B_R\setminus K$ and with boundary data $f$. To verify this claim, using the monotonicity of $\lambda$ over time, it is enough to check that  
\begin{equation}%\label{weak_limit}
\int_D  -\nabla p_0 \nabla\psi + \lambda(\cdot,0) \rho_0 \psi dx \geq 0
\end{equation}
for any nonnegative test function $\psi\in C^{\infty}_0(D)$. Since $\rho_0=\chi_{\{p_0>0\}}$, the question boils down to the nonnegativity of the measure $\mu_0:=\Delta p_0  + \lambda(\cdot,0)\chi_{\{p_0>0\}}.$ This follows the same proof of showing $\mu_t \geq 0$ in Proposition~\ref{prop:mu}, see section 5.

\medskip

 With the claim and the comparison principle for \eqref{eq:P} (Proposition~\ref{prop:order}), it follows that 
\begin{equation}\label{initiale} 
\rho(x,0) \leq \rho(x,\e) \hbox{ for all }\e>0.
\end{equation}
  Note that, since $\lambda$ is non-decreasing in time,  $\rho(\cdot,t-\e)$ is a subsolution of \eqref{eq:P}  for any  $\e>0$.  Thus by comparison principle and \eqref{initiale} it follows that 
$$
\rho(x, t-\e) \leq \rho(x,t) \hbox{ for any } t>\e>0,
$$
and we conclude that $\rho$ increases for all times.  $p$ accordingly increases by its definition.

\end{proof}

\begin{comment}

 For general  $\Sigma$ we use $L^1$ contraction.  Suppose $(\rho_1,p_1)$ and $(\rho_2,p_2)$ are two weak solutions of \eqref{eq:P} in  $\Omega \times \RR^+$. From Lemma~\ref{lem:basic} it follows that $\rho_i,p_i's$ are compactly supported, and $\rho_i(\cdot,t) \equiv 1$ in a small neighborhood $U$ of $K$. Thus from Proposition~\ref{prop:obstacle} it follows that  $p_i(\cdot,t)'s$ are smooth in $U$ for all $t>0$.  Note that, if $\rho_1 \geq \rho_2$, then $\partial_{\nu} p_1 \leq \partial_{\nu} p_2$ on $\partial K$ where $\nu$ is the inward normal,  since $p_1=p_2$ on $K$.  Therefore, from the weak equation, we have
$$
\begin{array}{lll}
\frac{d}{dt}[\int_{\RR^n} (\rho_1 - \rho_2) dx ] &=& \int_{\partial K} \partial_{\nu} p_1 - \partial_{\nu} p_2   dS - \lambda \int_{\RR^n}(\rho_1 - \rho_2) dx\\ \\
& \leq & - \lambda \int_{\RR^n} (\rho_1 - \rho_2) dx.
\end{array}
$$

Therefore it follows that, if $\Sigma$ can be approximated with a sequence of sets $\Sigma_n$ with $C^2$ boundaries, the statement is true. \textcolor{black}{how general $\Sigma$ can be? Certainly it includes $\Sigma$  with finite perimeter.}
\end{comment}

\begin{corollary}
Let $(\rho,p)$  be the weak solution of \eqref{eq:P} in $\Omega\times [0,\infty)$ with the fixed boundary data $p=f >0$ and the initial data $\rho_0\in BV$.  Then 
 $\Sigma(t):=\{p(\cdot,t)>0\}$ increases in time, and is a set of finite perimeter for a.e. $t>0$. Moreover for all $t\geq 0$
 \begin{equation}\label{presentation}
 \rho(\cdot,t) = \chi_{\Sigma(t)} + \rho^E\chi_{\RR^n\setminus \Sigma(t)}, \hbox{ where } \rho^E(x,t):= \rho_0\exp^{\int_0^t \lambda(x,s) ds}.
 \end{equation} 
\end{corollary}

\begin{proof}

We claim that the pressure support $\Sigma(t):=\{p(\cdot,t)>0\}$ increases over time.  For any $t_0> 0$, Let us call $\rho^*$ be the weak solution of \eqref{eq:P} with the initial data  $\chi_{\Sigma(t_0)}$, and with the same fixed boundary data $f$ for the pressure.  Then $\rho^*$ increases in time due to Lemma~\ref{patch}. From the monotonicity of $\rho^*$ and Proposition~\ref{prop:order}, we have 
\begin{equation} \label{order}
\chi_{\Sigma(t_0)} \leq \rho^*(\cdot,t) \leq \rho(\cdot,t_0+t) \hbox{ for all } t>0.
\end{equation}

It follows that $\Sigma(t)$ increases over time.  It follows from Proposition 1.5 that $\rho_t =  \lambda \rho$ in $\overline{\Sigma(t)}^C \times [0,t]$ for any $T>0$, and thus we can conclude \eqref{presentation}. Lastly $\Sigma(t)$, is a set of finite perimeter for a.e. $t>0$ since $\rho \in BV(\Omega)$ for a.e. $t>0$ and $\rho$ has jump discontinuity on the boundary of $\Sigma(t)$ due to \eqref{presentation}. 

\end{proof}

\appendix

\section{Tumor growth model with Nutrient}\label{app:tumor}

In \cite{PQV} (see also \cite{DP}), the following model for tumor growth is studied:
\begin{equation}
\begin{cases}
\pa_t \rho_m - \div( \rho_m\na p_m) = \rho_m G(p_m,c_m) \qquad x\in \RR^n, \; t\geq 0\\
\pa_t c_m -\Delta c_m + \rho_m H(c_m) = (c_B -c_m) K(p_m) \\
c_m(x,t)\to c_B \mbox{ for } x\to\infty
\end{cases}
\end{equation}
where
$$ p_m = \frac {m}{m-1} \rho_m^{m-1}.$$
In this system, the evolution of the cell population density $\rho_m\geq 0$ is coupled to the concentration of nutrients $c_m\geq 0$ by the cell division rate  $G(p,c)$. Importantly, this function satisfies
$$ \pa_ p G <-\beta<0$$
(see \cite{PQV} for a complete list of the assumptions necessary to get a good existence and uniqueness framework as well as the appropriate estimates to pass to the limit).

It is proved in \cite{PQV} that $\rho_m(x,t)$, $p_m(x,t)$ and $c_m(x,t)$ converge strongly in $L^1(Q_T)$ (for all $T>0$) to  $\rho_\infty,p_\infty, c_\infty$ in $BV(Q_T)$which  solves the system
\begin{equation}\label{eq:tumor}
\begin{cases}
\pa_t \rho_\infty - \div( \rho_\infty\na p_\infty) = \rho_\infty G(p_\infty,c_\infty) \qquad x\in \RR^n, \; t\geq 0\\
\pa_t c_\infty -\Delta c_\infty + \rho_\infty H(c_\infty) = (c_B -c_\infty) K(p_\infty) \\
c_\infty(x,t)\to c_B \mbox{ for } x\to\infty
\end{cases}
\end{equation}
with the Hele-Shaw relation  $p_\infty \in P_\infty(\rho_\infty)$.

Remarkably, the solution of this system is unique, and one would like to interpret the system as a weak form of some geometric Hele-Shaw type free boundary problem. For this one needs to identify the pressure $p_\infty$ as solution of an elliptic equation in $\{\rho_\infty=1\}$.

In \cite{DP}, it is proved that $p_\infty$ solves the complementarity condition 
$$ p_\infty (\Delta p_\infty + G(p_\infty,c_\infty))=0 \quad \mbox{ in } \mathcal D'(Q).$$ 
This condition says that $p_\infty$ solves an elliptic equation in $\{p_\infty\}$ and is proved by deriving additional estimates on $p_m$.

We will show below that the approach used in this paper can be used to characterize $p_\infty(\cdot,t)$ as the unique solution of an obstacle problem.
First, we summarize the estimates proved in \cite{PQV}:
\begin{lemma}
Under the assumptions listed in \cite{PQV}, the following holds for all $T>0$:
\begin{itemize}
\item $\rho_m(t)$ is uniformly compactly supported for $t\in [0,T]$;
\item $|\na p_m|$ is bounded in $L^2(Q_T)$
\item $0\leq p_m\leq p_M$, $0\leq \rho_m\leq \left(\frac{m-1}{m} p_M\right)^\frac{1}{m-1}$, $0<c_m<c_B$
\item $\rho_m$, $p_m$ and $c_B-c_m$ are bounded in $BV(Q_T)$
\item $\rho_m$, $p_m$ and $c_B-c_m$ converge   strongly in $L^1$ and almost everywhere to $\rho_\infty$, $p_\infty$ and $c_B-c_\infty$.
\end{itemize}
\end{lemma}
Furthermore, proceeding as in Lemma \ref{lem:2.4}, it is not difficult to show that $\{\rho_m\}_{m\in \NN}$ is relatively compact in $C^s(0,T;H^{-1}(\RR^n))$ for all $s\in(0,1/2)$ and thus that $\rho_\infty \in C(0,T;H^{-1}(\RR^n))$. 

Finally, since $p_\infty$ and $c_B-c_\infty$ are in $BV(Q_T)$, we can define the trace $p^+(\cdot,t)$ and $c^+(\cdot,t)$ for all $t>0$ as in \eqref{eq:tracep}.
We can then prove the following result:
\begin{proposition}\label{prop:obstumor}
For all $t>0$, let  $E_t$ denote the space
$$ E_t = \{ v\in H^1(\RR^n)\cap L^1(\RR^n)\, ;\,  \; v (x) \geq 0, \;    \langle v, 1-\rho_\infty(t) \rangle_{H^1,H^{-1}} =0 \}.$$
Then {for all $t>0$}, the function $x\mapsto p^+(x,t)$ is the unique solution of the minimization problem:
\begin{equation}\label{eq:obstacletumor}
\begin{cases}
p \in E_t \\ 
\displaystyle \int_{\RR^n} \frac 1 2 |\na p|^2 -\G(p,c^+) \, dx \leq \int_{\RR^n} \frac 1 2 |\na v|^2 -\G( v,c^+)\, dx \qquad \forall v \in E_t
\end{cases}
\end{equation}
where $\G$ is the (concave) function such that $\pa_p \G(p,c) = G(p,c)$ and $\G(0,c)=0$.
Furthermore $p_\infty$ satisfies the complementarity condition
\begin{equation}\label{eq:comp}
p_\infty (\Delta p_\infty + G(p_\infty,c_\infty) ) =0 \quad \mbox{ in } \mathcal D'(\RR^n\times(0,\infty)).
\end{equation}
\end{proposition}
As mentioned in the introduction (see Proposition \ref{prop:comp}), if the complementarity condition \eqref{eq:comp}  is known to hold, then one can derive the variational formulation  \eqref{eq:obstacletumor} from the weak equation \eqref{eq:tumor}.
In particular, this   complementarity condition was derived for this particular model in \cite{DP} by using a generalized Aronson-B\'enilan estimate and the $L^2(W^{1,4})$ estimate on the pressure (but our proof here does not require either of these estimates). 
\begin{proof}
First we recall the equation for the pressure $p_m$:
\begin{equation}\label{eq:pressuretumor}
\pa_t p_m = (m-1) p_m( \Delta p_m + G(p_m,c_m)) + |\na p_m|^2.
\end{equation}
We then proceed as in  the proof of Theorem \ref{prop:obs}:
Given $t_0>0$ and  a function $v(x)$ in $E_{t_0}$, we use the equation for the pressure \eqref{eq:pressuretumor} and density \eqref{eq:tumor} to write:
\begin{align*}
&  \int_{\RR^n} \na p_m \cdot \na p_m - \rho_m \na p_m\cdot\na v- \G(p_m,c_m) + \G(v,c_m)\, dx\\
& \qquad \qquad = -\frac{1}{m-1}\left[ \frac{d}{dt}\int_{\RR^n} p_m\, dx - \int_\Omega|\na p_m|^2\, dx\right]\\
&  \qquad \qquad \quad + \frac{d}{dt} \int_{\RR^n} v \rho_m\, dx\\
&  \qquad \qquad\quad + \int_{\RR^n} p_m G(p_m,c_m) - \rho_m v G(p_m,c_m) - \G(p_m,c_m) + \G(v,c_m)\, dx
\end{align*}
in $\mathcal D'(\RR_+)$. 
Using the concavity of $\G$ to write 
$$ \G(v,c_m) - \G(p_m,c_m) \leq G(p_m,c_m)(v-p_m)$$
we deduce
\begin{align*}
  \int_{\RR^n} \na p_m \cdot \na p_m - \rho_m \na p_m\cdot\na v- \G(p_m,c_m) + \G(v,c_m)\, dx& = -\frac{1}{m-1}\left[ \frac{d}{dt}\int_{\RR^n} p_m\, dx - \int_{\RR^n}|\na p_m|^2\, dx\right]\\
&  \quad + \frac{d}{dt} \int_{\RR^n} v \, \rho_m\, dx + \int_{\RR^n} (1- \rho_m) \, v \, G(p_m,c_m) \, dx.
\end{align*}

We can now proceed as in the proof of Theorem \ref{prop:obs}:
Integrating this equality with respect to $t\in [t_0,t_0+\delta)$ and 
using the weak $L^2$ convergence of $\na p_m$ and $\rho_m \na p_m$ to $\na p$, we get 
\begin{align*}
& \int_{t_0}^{t_0+\delta} \int_{\RR^n} |\na p_\infty|^2 -  \na p_\infty \cdot \na v- \G(p_\infty,c_\infty) + \G(v,c_\infty)\, dx\, dt\\
& \quad \leq \int_{\RR^n} v(x) [\rho_\infty(x,t_0+\delta) - \rho_\infty(x,t_0)] \, dx+ \int_{t_0}^{t_0+\delta} \int_{\RR^n}  (1- \rho_\infty) v G(p_\infty,c_\infty) \, dx\, dt\\
& \quad \leq \| G(p_\infty,c_\infty)\|_{L^\infty}  \int_{t_0}^{t_0+\delta} \int_{\RR^n}  v(1- \rho_\infty )   \, dx\, dt
%& \quad \leq  \int_{t_0}^{t_0+\delta} \int_{\RR^n}  (v- \rho_\infty v)  G(p_\infty,c_\infty) \, dx\, dt
\end{align*}
(where we used the fact that $v(x) \rho_\infty(x,t_0) = v(x)$ and $v(x) \rho_\infty(x,t) \leq v(x)$ for all $t$)

Finally, dividing by $\delta$  and using Young's inequality, we rewrite the inequality as
\begin{align*}
& \frac 1 \delta  \int_{t_0}^{t_0+\delta} \int_{\RR^n}\frac 1 2  |\na p_\infty|^2 - \G(p_\infty,c_\infty) \, dx\\
&\qquad \leq 
   \frac 1 \delta \int_{t_0}^{t_0+\delta} \int_{\RR^n}  \frac 1 2  |\na v|^2 -  \G(v,c_\infty)\, dx\, dt 
   +   \frac C \delta\int_{t_0}^{t_0+\delta} \int_{\RR^n}  v (1- \rho_\infty )   \, dx\, dt
   \\
 & \qquad \leq  \int_{\RR^n}  \frac 1 2  |\na v|^2 - \frac 1 \delta\int_{t_0}^{t_0+\delta}  \G(v,c_\infty)\, dt\, dx
+   \frac C \delta \int_{t_0}^{t_0+\delta} \langle v, 1- \rho_\infty \rangle_{H^1,H^{-1}}    \, dx\, dt.
\end{align*}
The continuity of $t\mapsto  \langle v, 1- \rho_\infty \rangle_{H^1,H^{-1}}$ and the fact that $v\in E_t$ implies that the last term converges to zero as $\delta\to0$.
{We can now conclude as in the proof of Theorem \ref{prop:obs}}.
%Lebesgue differentiation theorem implies that the left hand side converges to
%$$\int_{\RR^n}\frac 1 2  |\na p_\infty(x,t_0) |^2  - \G(p_\infty(x,t_0),c_\infty(x,t_0))  \, dx$$
%for a.e. $t_0>0$ and the result follows.
\medskip

Finally, given a test function $\vphi\in \mathcal D(\RR^n\times(0,\infty))$, we take $v = p_\infty + \eps (p_\infty\vphi) =p_\infty(1+\eps\vphi)$  in \eqref{eq:obstacletumor}, with $|\eps|$ small enough so that $1+\eps\vphi \geq0$.
Passing to the limit $\eps\to 0^-$ and $\eps\to 0^+$ yields
$$
 \int_{\RR^n}  \na p_\infty\cdot \na (p_\infty\vphi) -G(p_\infty,c_\infty) p_\infty\vphi\, dx =0$$
and \eqref{eq:comp} follows.

\end{proof}

\section{The complementarity condition}

\begin{proof}[Proof of Proposition \ref{prop:comp}]
We note that $ \pa_t \rho = \Delta p + \lambda \rho \in L^2(0,T;H^{-1} (\Omega))$.
Given $u\in E_t$, we have $p-u \in L^2(0,T;H^1_0(\Omega))$ and so we can write (in $\mathcal D'(\RR_+)$):
\begin{equation}\label{eq:fe}
\langle \pa_t \rho, (p-u)\rangle_{H^{-1}, H^1_0} =  \langle \Delta p + \lambda \rho, p-u\rangle_{H^{-1}, H^1_0}= - \int_\Omega \na p \cdot \na (p-u) - \lambda \rho (p-u)\, dx .
\end{equation}

Next, proceeding as in the beginning of the proof of Lemma \ref{patch} (using the  comparison principle for the limiting problem, Proposition~\ref{prop:order}), we can show that $\rho=1$ in $U\times \RR_+$ for some neighborhood $U$ of $K$ and that $\supp p$ is bounded in $\Omega \times [0,T]$.
In particular,  $\pa_t \rho$ vanishes in $U\times \RR_+$.
Taking a smooth function $\phi(x)$ which is equal to $1$ in $\supp p\setminus (U\times [0,T])$ and vanishes on $\pa K$, we can write
\begin{align*}
\langle \pa_t \rho, (p-u)\rangle_{H^{-1}, H^1_0} & = \langle \pa_t \rho, (p-u)\phi \rangle_{H^{-1}, H^1_0}\\
& =  \langle \Delta p + \lambda \rho , p \phi \rangle_{H^{-1} , H^1_0 } -\langle \pa_t \rho, u\phi \rangle_{H^{-1}, H^1_0} \\
& =  \langle p( \Delta p + \lambda \rho)  , \phi \rangle_{\mathcal D',\mathcal D}  -\langle \pa_t \rho, u\phi \rangle_{H^{-1}, H^1_0} \\
& =  \langle p( \Delta p + \lambda \rho)  , \phi \rangle_{\mathcal D',\mathcal D}- \frac{ d}{dt} \int_\Omega \rho u \phi \, dx\\
& = - \frac{ d}{dt} \int_\Omega \rho u \phi \, dx
\qquad \mbox{ in $\mathcal D'(\RR_+)$}
\end{align*}
where we used the fact that $ \langle p( \Delta p + \lambda \rho)  , \phi \rangle_{\mathcal D',\mathcal D}=0$ (this is the complementarity condition). Using \eqref{eq:fe}, we deduce
$$  \int_\Omega \na p \cdot \na (p-u) - \lambda \rho (p-u)\, dx= 
 \frac{ d}{dt} \int_\Omega \rho u \phi \, dx
\qquad \mbox{ in $\mathcal D'(\RR_+)$.}
$$
Using the fact that $\rho(x,t)p(x,t) = p(x,t)$, we deduce:
\begin{align*}
\int_\Omega  \na p \cdot \na (p-u)- \lambda  (p-u)\, dx \, dt 
& =
 \int_\Omega  \na p \cdot \na (p-u) - \lambda \rho (p-u)\, dx \, dt  +
    \int_\Omega \lambda (1-\rho) u\, dx \, dt\\
& \leq 
  \frac{ d}{dt} \int_\Omega \rho u \phi \, dx +
\Lambda  \int_\Omega  (1-\rho) u\, dx \, dt\\
\end{align*}
Integrating with respect to $t\in [t_0,t_o+\delta]$, we get
\begin{align*}
\int_{t_0}^{t_0+\delta} 
\int_\Omega \na p \cdot \na (p-u) - \lambda  (p-u)\, dx \, dt 
%& \leq 
%\int_{t_0}^{t_0+\delta} 
%\int_\Omega \na p \cdot \na (p-u) - \lambda \rho (p-u)\, dx \, dt \\
%&\qquad +\Lambda \int_{t_0}^{t_0+\delta}  \int_\Omega \lambda (1-\rho) u\, dx \, dt\\
& \leq 
 \int_\Omega (\rho(t_0+\delta)-\rho(t_0)) u\phi\, dx+ \Lambda \int_{t_0}^{t_0+\delta}  \int_\Omega  (1-\rho) u\, dx \, dt\\
& \leq 
 \int_\Omega (\rho(t_0+\delta)-1) u\phi\, dx+\Lambda \int_{t_0}^{t_0+\delta}  \int_\Omega  (1-\rho) u\, dx \, dt\\
& \leq \Lambda \int_{t_0}^{t_0+\delta}  \int_\Omega  (1-\rho) u\, dx \, dt
\end{align*}
and the result now follows by proceeding as in the proof of Theorem \ref{prop:obs}.
\end{proof}

\bibliographystyle{alpha}
\bibliography{bib_PME}

\end{document}